\documentclass[review, 11pt]{elsarticle}

\usepackage[utf8]{inputenc}
\usepackage[top=2cm, bottom=2cm, left=2cm, right=2cm]{geometry}
\usepackage{float}
\usepackage{caption}
\usepackage{tikz,pgfplots,pgfplotstable,booktabs}
\usepackage[official]{eurosym}
\usepackage{amsthm}
\usepackage{xcolor}
\usepackage{bm}
\usepackage{colortbl}
\pgfplotsset{compat=1.15}
\usepackage{comment}
\usepackage{circuitikz}
\usepackage{tikz,pgfplots,p   gfplotstable,booktabs}%
\usepackage{multirow, array}
\usepackage{eurosym}
\usepackage[mathscr]{euscript}
\usepackage{amsmath,amssymb,amsfonts}
\usepackage{algorithmic}
\usepackage[Algorithm]{algorithm}
\usepackage{graphicx}
\usepackage{textcomp}
\usepackage{appendix}
\usepackage{mathtools}
\DeclarePairedDelimiter\ceil{\lceil}{\rceil}
\DeclarePairedDelimiter\floor{\lfloor}{\rfloor}
\usepackage{hyperref}
\usepackage{mathabx}
\usepackage[subrefformat=parens,labelformat=parens]{subfig}
\hypersetup{
    colorlinks=false,
    hidelinks
}

\newtheorem{theorem}{Theorem}[section]
\newtheorem{proposition}[theorem]{Proposition}
\newtheorem{corollary}[theorem]{Corollary}

\newtheorem{definition}[theorem]{Definition}

\allowdisplaybreaks

\allowdisplaybreaks

\makeatletter
\newcommand*\bigcdot{\mathpalette\bigcdot@{.5}}
\newcommand*\bigcdot@[2]{\mathbin{\vcenter{\hbox{\scalebox{#2}{$\m@th#1\bullet$}}}}}
\makeatother

\usepackage{csquotes}

\usepackage{natbib}
 \bibpunct[, ]{(}{)}{,}{a}{}{,}%
\journal{arXiv}

\begin{document}

\begin{frontmatter}

\title{Tight and Compact Sample Average Approximation for Joint Chance-constrained Problems with Applications to Optimal Power Flow}


\author[mymainaddress]{Álvaro Porras}
\author[mymainaddress]{Concepción Domínguez}
\author[mymainaddress]{Juan Miguel Morales\corref{mycorrespondingauthor}}
\cortext[mycorrespondingauthor]{Corresponding author}
\ead{juan.morales@uma.es}
\author[mymainaddress]{Salvador Pineda}

\address[mymainaddress]{Optimization and Analytics for Sustainable energY Systems (OASYS), University of Málaga, 29071 Málaga (Spain)}

\begin{abstract}
In this paper, we tackle  the resolution of chance-constrained problems reformulated via Sample Average Approximation. The resulting data-driven deterministic reformulation takes the form of a large-scale mixed-integer program cursed with Big-Ms. We  introduce an exact resolution method for the MIP that combines the addition of a set of valid inequalities to tighten the linear relaxation bound with coefficient strengthening and constraint screening algorithms to improve its Big-Ms and considerably reduce its size. The proposed valid inequalities are based on the notion of \emph{k-envelopes}, can be computed offline using polynomial-time algorithms, and added to the MIP program all at once. Furthermore, they are equally useful to boost the strengthening of the Big-Ms and the screening rate of superfluous constraints. We apply our procedures to a probabilistically-constrained version of the DC Optimal Power Flow problem with uncertain demand. The chance constraint requires that the probability of violating \emph{any} of the power system's constraints be lower than some parameter $\epsilon > 0$. In a series of numerical experiments which involve five power systems of different size, we show the efficiency of the proposed methodology and compare it with some of the best-performing convex inner approximations currently available in the literature. 
\end{abstract}

\end{frontmatter}

\section{Introduction} 

Chance-constrained programming suits applications in areas where decisions have to be made dealing with random parameters (\cite{miller1965chance}). In these situations, it is desirable to ensure feasibility of the system almost surely, but there is hardly any decision which would guarantee it under extreme events or unexpected random circumstances. Within this context, chance-constrained programming offers a powerful modeling framework to identify cost-efficient decisions that satisfy the problem's constraints with a high level of confidence. 

In this paper, we focus on chance-constrained problems (CCPs) with joint linear chance constraints and random RHS and LHS, i.e.,
\begin{subequations} \label{GenCC}
\begin{align}
 \min_{x} \quad &  f(x) \label{GenCC_FO}\\
\text{s.t.}  \quad & x \in X \label{GenCC_xinX} \\
& \mathbb{P}\left\{ a_{j}(\omega)^{\top}x \le b_j(\omega), \ \forall j \in \mathcal{J}\right\} \ge 1-\epsilon. \label{GenCC_chance}
\end{align}
\end{subequations}

In \eqref{GenCC}, $x\in \mathbb{R}^{|\mathcal{I}|}$ is a vector of continuous decision variables, $X\subseteq \mathbb{R}^{|\mathcal{I}|}$ represents a set of deterministic constraints, and $f: \mathbb{R}^{|\mathcal{I}|} \longrightarrow \mathbb{R}$ is a convex function. Uncertainty is represented through the random vector $\omega$ taking values in $\mathbb{R}^{d}$ and giving rise to a technology matrix with random rows $a_{j}(\omega)\in \mathbb{R}^{|\mathcal{I}|}$, $j \in \mathcal{J}$ and random $b_j(\omega) \in \mathbb{R}$, $j \in \mathcal{J}$. 
$\mathbb{P}$ is a probability measure, and $\epsilon$ is a confidence or risk parameter, typically near zero, so that the set of constraints \eqref{GenCC_chance} are satisfied with probability at least $(1-\epsilon)$.  Apart from power systems, applications of CCPs include supply chain, location and logistics (\cite{taleizadeh2012, shaw2016, elci2018a}), risk control in finance (\cite{danielsson2008, natarajan2008}), and healthcare problems such as operating room planning (\cite{najjarbashi2020}) or vaccine allocation (\cite{tanner2008}), among others.

When the probabilistic constraint corresponds to \eqref{GenCC_chance}, the CCP has joint chance constraints (JCC) and is hence classified as a joint CCP (JCCP), in contrast with single CCPs (SCCPs), i.e.\ CCPs with individual or single chance constraints (SCC) of the form $\mathbb{P}\left\{ a_{j}(\omega)^{\top}x - b_j(\omega) \le 0\right\} \ge 1-\epsilon_j$, $\forall j \in \mathcal{J}$. JCCPs are suitable for contexts where all constraints need to be simultaneously satisfied with a high probability, and the dependence between random variables makes them clearly harder. Both SCCPs and JCCPs have been extensively studied (see \cite{prekopa2003, vanackooij2011} and the references therein).

There are a number of reasons why general CCPs are challenging. The first one is the non-convexity of the feasible set. In general, the feasible region of a CCP is not convex in the original space even when $x$ is continuous, there is only RHS uncertainty and the constraints inside the probability in \eqref{GenCC_chance} define a polyhedral region (\cite{kucukyavuz2022}). To circumvent this problem, several approaches have been proposed. Some methods (e.g.\ \cite{lagoa2005,henrion2007,henrion2011}) give convexity results and investigate the conditions under which the feasible region of problem \eqref{GenCC} is convex. In another line of research, various convex approximation schemes such as quadratic (\cite{ben-tal2000}) or Bernstein approximation (\cite{nemirovski2007}), have been proposed in the literature. The CVaR approximation has gained a lot of popularity since its introduction (\cite{rockafellar2000, sun2014}), and remains one of the most used methods to deal with stochastic problems.  Nonetheless, the solutions to the approximated problems err on the side of over-conservatism. In this context, some iterative schemes such as ALSO-X have been recently proposed to identify tighter inner convex approximations of the CCP at the expense of a higher computational cost (\cite{ahmed2017, jiang2022}). Finally, other works suggest convex approximations for non-linear CCPs. For instance, \cite{hong2011} propose to solve the JCCP by a sequence of convex approximations followed by a gradient-based Monte Carlo method, whereas \cite{pena-ordieres2020} introduce a smooth sampling-based approximation.

The second difficulty of CCPs is that checking the feasibility of a given solution is not, in general,  an easy task. For instance, even if the uncertainty follows a known continuous distribution, calculating the joint probability requires a multi-dimensional integration, which becomes increasingly difficult with the dimension of the random vector~$\omega$. On top of that, in most cases the distribution  $\mathbb{P}$ is not fully known.

A commonly used approach to cope with these obstacles is to replace $\mathbb{P}$ with its empirical estimate based on a set of random samples (\cite{shapiro2003}). This approach is known in the technical literature as Sample Average Approximation (SAA) and may also be seen as a data-driven approach that works with observations of $\omega$ that are available to the decision-maker even if the true data-generating distribution $\mathbb{P}$ is unknown. The main advantage of  SAA  is that it allows to reformulate CCPs as deterministic optimization problems that can be solved using standard optimization techniques.

Indeed, the deterministic reformulation of  \eqref{GenCC} using SAA results in a mixed-integer problem (MIP). For the resolution of CCPs using MIP reformulations, we refer the reader to \cite{kucukyavuz2022}. When there is only RHS uncertainty (i.e.\ the technology matrix is fixed), a reformulation of the problem leads to a MIP with a set of constraints that form a \emph{mixing set} and that has been extensively studied, alone or in combination with the knapsack constraint that also appears in the formulation (\cite{gunluk2001,luedtke2010, abdi2016}). Alternative reformulations like the ones proposed by \cite{dentcheva2000, nair2011} rely on the concept of $(1-\epsilon)$-efficient points. The case when the technology matrix is random, while the RHS is not, has been studied e.g.\ in \cite{tayur1995}. As for the general case, it has also been addressed in the literature. Specifically, a large line of research has focused on the development of \emph{quantile cuts}, a particular type of valid inequality that can be viewed as a projection of a set of mixing inequalities for the MIP onto the original problem space. These cuts and the associated quantile closure have been recently studied in \cite{qiu2014,xie2016,xie2018,ahmed2018} and successfully applied to computational experiments of CCPs in \cite{song2014,ahmed2017}, among others.

Regarding its disadvantages, the SAA solution tends to be optimistic and underestimate the probability of violating the constraints (\cite{nemirovski2006}). In fact, the performance of SAA is directly contingent on the number of samples available. We refer the reader to the article by \cite{luedtke2008}, where relations are established between the empirical acceptable probability of violation and the number of samples such that the SAA-based solution be feasible in the JCCP with a predefined confidence level. In this paper, however, we do not investigate how well the SAA-based solution approximates the solution of the original JCCP. Statistical considerations apart, the main focus of this paper is on proposing an efficient methodology to solve the deterministic optimization model that results from reformulating JCCPs using SAA. In particular, the method we propose solves the MIP to optimality, and is based on the combination of a tightening-and-screening procedure with the development of valid inequalities to obtain a formulation which is compact and tight at a time. 

We begin with a description of an iterative algorithm to strengthen the Big-Ms present in the mixed-integer reformulation of the JCC. Interestingly, although the procedure is not new (see \cite{qiu2014}), we complement it with a screening procedure that allows us to eliminate inequalities of the MIP. 

Next, we introduce a new set of valid inequalities designed to strengthen the linear relaxation of the model. As for their structure, each valid inequality is developed using quantile information and avoiding the use of Big-M constants, which are known to lead to weak linear relaxation bounds. Furthermore, from a theoretical point of view, they are closely related to the aforementioned quantile cuts. In fact, as proven in \ref{sec:anexoQuantileCuts}, the addition of our inequalities yields a feasible region equivalent to the quantile closure of a specific relaxation of the MIP problem. The main advantage of our inequalities is that, unlike quantile cuts, they are not NP-hard to compute,  and neither do they require a specific separation algorithm to include them dynamically (that is, they can all be included in the model from the outset). As many other techniques developed for JCCPs, our cuts also apply to the SCCPs. Finally, our valid inequalities extend the results introduced in \cite{roos1994} for the $k$-\emph{violation problem}, which can be seen as a mixed-integer linear problem (MILP) reformulation of a SCCP. For ease of explanation, the valid inequalities are proposed for a particular class of linear chance constraints, but our results are also applicable to more general types of CCPs, as we detail in \ref{sec:anexoExtensionMultidim}. 

We apply our resolution method to the Optimal Power Flow (OPF) problem, for which we propose a chance-constrained formulation and extensive computational experiments featuring results for five standard power systems available in the literature. The combination of the valid inequalities with the tightening of the Big-Ms and the screening procedure allows us to effectively solve to optimality instances that are not solved with the initial MIP formulation, since the combination of both techniques notably reduces their size and difficulty. We also compare our resolution approach with state-of-the-art convex inner approximations of CCPs, in particular, the CVaR-based approximation, ALSO-X, and ALSO-X+ (\cite{jiang2022}). Finally, we compare our results with an exact approach for solving CCPs assuming SAA, namely, the branch-and-cut decomposition algorithm proposed in \cite{luedtke2014}. This algorithm is based on the addition of quantile cuts strengthened using the star inequalities of \cite{atamturk2000} to a relaxation of the original problem.

The remainder of the paper is organized as follows. Section~\ref{sec:JCCP} involves the reformulation of the problem into a MIP using  SAA  and the proposed methodology: Subsection~\ref{sec:screening} describes the tightening and screening procedures, whereas Subsection~\ref{sec:Valid} introduces the valid inequalities and the necessary algorithms to compute them. Section~\ref{sec:OPF} is devoted to the application of the results to the OPF problem. The section begins with an introduction to the  problem and the associated literature review, and Subsection~\ref{sec:OPF_notation} formalizes the notation and the mathematical model solved in the computational study. In Subsection~\ref{sec:OPF_casestudy} we present the case study, testing our results to solve instances of the DC OPF available in the literature. Section~\ref{sec:conclusion} points further research topics and includes some concluding remarks. Finally, the main text is accompanied by two appendices: \ref{sec:anexoQuantileCuts} establishes the relationship between our valid inequalities and the quantile cuts present in the literature, and \ref{sec:anexoExtensionMultidim} details the complexity of the algorithms proposed in Subsection \ref{sec:Valid} and discusses the generalization of our methodology and its possible application to other types of CCPs.

\section{Solving JCCPs via Sample Average Approximation} \label{sec:JCCP}
Using SAA, problem \eqref{GenCC} can be reformulated into a MIP. To this aim, assume that $\omega$ has a finite discrete support defined by a collection of points $\{\omega_s \in \mathbb{R}^{|\mathcal{N}|}, s \in \mathcal{S}\}$ and respective probability masses $\mathbb{P}(\omega = \omega_s)=\frac{1}{|\mathcal{S}|}$, $\forall s\in \mathcal{S}=\{1,\dots,|\mathcal{S}|\}$. Consider the rows $a_{js}=a_j(\omega_s)$, and the random $b_{js}=b_j(\omega_s)$. Define $p:=  \floor*{\epsilon |\mathcal{S}|}$, the vector $y$ of binary variables $y_s$, $\forall s\in \mathcal{S}$, and the large enough constants $M_{js}$. Then, the MIP reformulation of problem \eqref{GenCC} can be stated as:
\begin{subequations}\label{MIP1GenCC}
\begin{align}
 \min_{x} \quad &  f(x) \label{MIP1GenCC_FO}\\
\text{s.t.}  \quad & x \in X  \label{MIP1GenCC_xinX}\\
&  a_{js}^{\top}x \le b_{js} + M_{js}y_s, \quad \forall j\in \mathcal{J}, s \in \mathcal{S}\label{MIP1GenCC_chance}\\
& \sum_{s \in \mathcal{S}} y_s \leq p \label{MIP1GenCC_violation}\\
&  y_s \in \{0,1\}, \quad \forall s \in \mathcal{S}. \label{MIP1GenCC_y}
\end{align}
\end{subequations}
Using the generic formulation \eqref{MIP1GenCC}, we present in Subsection~\ref{sec:screening} a procedure to properly tune the values of the large constants $M_{js}$. We also explain in this subsection how the intermediate results of the tightening procedure can be efficiently used to remove constraints from set \eqref{MIP1GenCC_chance} that are superfluous, thus making model \eqref{MIP1GenCC} more compact. Finally, we introduce in Subsection~\ref{sec:Valid} a set of valid inequalities that makes the linear relaxation of a particular case of~\eqref{MIP1GenCC} remarkably tighter, and we relate our inequalities to the \textit{quantile closure} studied in \cite{xie2018}. 
\subsection{Tightening and Screening}
\label{sec:screening}
It is well-known that the linear relaxation of a Big-M formulation tends to provide weak lower bounds in general (\cite{conforti2014}). This is even more so when the Big-Ms are chosen too loose. Constants $M_{js}$ in \eqref{MIP1GenCC_chance} should be set large enough for the corresponding constraints to be redundant when the associated binary variables $y_s$ are equal to 1, and \emph{as small as possible} to tighten the MIP formulation. To this end, \cite{qiu2014} provide an algorithm called ``Iterative Coefficient Strengthening'' that has been successfully applied to other CCPs (see, e.g.,\ \cite{song2014}). A customization of this procedure for the joint chance-constrained formulation \eqref{MIP1GenCC} is detailed in Algorithm~\ref{alg:IterativeCoefStreng}.

\begin{algorithm}
\begin{small}
\caption{Iterative Coefficient Strengthening} \label{alg:IterativeCoefStreng}
\begin{algorithmic} 
\STATE \textbf{Input}: The LHS and RHS vectors  $\{a_{js}\}_{j \in \mathcal{J}, s \in \mathcal{S}}$ and coefficients $\{b_{js}\}_{j \in \mathcal{J}, s \in \mathcal{S}}$, respectively, and the maximum allowed violation probability $p$ that determine the joint chance-constraint system \eqref{MIP1GenCC_violation}, the deterministic feasible set $X$, and the total number of iterations $\kappa$.

\textbf{Output}: The large constants $M_{js}$, $\forall j\in \mathcal{J}$, $s\in \mathcal{S}$.

\begin{enumerate}[\hspace*{0.5cm} Step 1.]\itemsep -.1cm
\item Initialization, $k = 0$, $M_{js}^0 = \infty$, $\forall j\in \mathcal{J}$, $s\in \mathcal{S}$.
\item For each $j\in\mathcal{J}$ and $s\in\mathcal{S}$ update $M_{js}^{k+1}$ as follows: If $M_{js}^{k}<0$, then $M_{js}^{k+1}=M^{k}_{js}$. Otherwise,
\begin{subequations}
\label{RelMIP}
\begin{flalign}
M^{k+1}_{js}=\max_{x,y_s} \quad &   a_{js}^{\top}x - b_{js} \label{RelMIP_FO}\\
\text{s.t.}  \quad & x \in X  \label{RelMIP_xinX}\\
&  a_{js}^{\top}x  - b_{js} \leq M^{k}_{js}y_s, \quad \forall j\in \mathcal{J}, s\in \mathcal{S} \label{RelMIP_chance}\\
& \sum_{s \in \mathcal{S}} y_s \leq p \label{RelMIP_violation}\\
& 0 \leq y_s \leq 1, \quad \forall s \in \mathcal{S}. \label{RelMIP_y}
\end{flalign}
\end{subequations}
\item If $k+1 < \kappa$, then $k = k + 1$ and go to Step 2. Otherwise, stop.
\end{enumerate}
\end{algorithmic}
\end{small}
\end{algorithm}

The output of Algorithm~\ref{alg:IterativeCoefStreng} is the tuned Big-Ms that are input to the MIP reformulation \eqref{MIP1GenCC}. This algorithm produces Big-Ms whose value either decreases or remains equal at each iteration. There is, in fact, a number of iterations beyond which the resulting Big-Ms converge. To reduce the computational burden of running Algorithm~\ref{alg:IterativeCoefStreng}, all problems \eqref{RelMIP} in Step 2 can be solved in parallel. Hereinafter, we use the short name ``\textbf{T}($\kappa$)'' (from ``Tightening'') to refer to the solution of model \eqref{MIP1GenCC} using the Big-M values given by the ``Iterative Coefficient Strengthening'' algorithm with $\kappa$ iterations.

Equally important, Algorithm~\ref{alg:IterativeCoefStreng} can be easily upgraded to delete constraints $(j,s)$ in \eqref{MIP1GenCC_chance} that are redundant, and therefore, can be removed from problem \eqref{MIP1GenCC}. Indeed, if Algorithm~\ref{alg:IterativeCoefStreng} delivers a large constant $M_{js} \leq 0$, then constraint $j$ in scenario $s$ can be deleted from \eqref{MIP1GenCC} without altering its feasible region or its optimal solution. This is so because a non-positive $M_{js}$ means that there is no $x$ satisfying \eqref{RelMIP_xinX}--\eqref{RelMIP_y} such that the constraint takes on a value strictly greater than zero. Consequently, the constraint is redundant in \eqref{MIP1GenCC}, since the feasibility region of \eqref{RelMIP} is a relaxation of \eqref{MIP1GenCC}. This upgrade of method \textbf{T} not only makes formulation \eqref{MIP1GenCC} tighter through coefficient strengthening, but also more compact by screening out redundant constraints. Naturally, the tightening and screening power of algorithm \textbf{T} increases at each iteration. From now on, we use the short name ``\textbf{TS}($\kappa$)'' (from ``Tightening and Screening'') to refer to the strategy whereby model \eqref{MIP1GenCC} is solved without the constraints \eqref{MIP1GenCC_chance} for which the value of $M_{js}$ provided by Algorithm~\ref{alg:IterativeCoefStreng} after $\kappa$ iterations is lower than or equal to 0.

While strengthening the parameters $M_{js}$ is a common strategy in the technical literature to reduce the computational burden of CCPs, this is the first time, to our knowledge, that intermediate results of Algorithm~\ref{alg:IterativeCoefStreng} are used to eliminate superfluous constraints from model \eqref{MIP1GenCC}. We stress that the screening process itself comes at no cost from Algorithm~\ref{alg:IterativeCoefStreng}, while removing superfluous constraints from \eqref{MIP1GenCC} may substantially facilitate its solution. As we show in Section~\ref{sec:OPF_casestudy}, this is particularly true for the JCC-OPF.

Apart from making formulation \eqref{MIP1GenCC} more compact, the screening of superfluous constraints can also be used to accelerate the ``Iterative Coefficient Strengthening'' process at each iteration. To do so, it suffices to modify Algorithm~\ref{alg:IterativeCoefStreng} so that \eqref{RelMIP_chance} only includes the constraints for which $M^k_{js}>0$. Thus, the number of constraints of model \eqref{RelMIP} is significantly reduced and so is its solution time. 

\subsection{Valid inequalities}
\label{sec:Valid}
In this subsection, we propose valid inequalities that apply to JCCPs with a particular structure: each row $a_{j}(\omega)$ in the technology matrix of \eqref{GenCC} can be rewritten as $a_j(\omega) = a^{0}_j + \Omega_j(\omega) \hat{a}_j$, with $a_j^{0}, \hat{a}_{j}\in \mathbb{R}^{|\mathcal{I}|}$ and where $\Omega_j(\omega_s)=\Omega_{js}$ is a real-valued function whose domain includes the support of $\omega$. In other words, we derive a set of valid inequalities to make the linear relaxation of problem
\begin{subequations}
\label{MIPGenCC}
\begin{align}
 \min_{x,y_s} \quad &  f(x) \label{MIPGenCC_FO}\\
\text{s.t.}  \quad & x \in X  \label{MIPGenCC_xinX}\\
& \Omega_{js} \hat{a}_j^{\top} x - b_{js} + x^{\top} a_j^{0} \le M_{js}y_s, \quad \forall j\in \mathcal{J}, s\in \mathcal{S} \label{MIPGenCC_chance}\\
& \sum_{s \in \mathcal{S}} y_s \leq p \label{MIPGenCC_violation}\\
& y_s \in \{0,1\}, \quad \forall s \in \mathcal{S} \label{MIPGenCC_ybinaria}
\end{align}
\end{subequations}
\noindent tighter. Note that we do not make any assumptions on the sign of $\hat{a}_j^{\top}$, $b_{js}$ and $a_j^{0}$. Furthermore, these valid inequalities can also be added to the constraint set~\eqref{RelMIP} of Algorithm~\ref{alg:IterativeCoefStreng}, thus dramatically increasing the tightening and screening power of \textbf{TS}. To facilitate the comparative analysis carried out in Section~\ref{sec:OPF_casestudy}, the so upgraded  algorithm is named ``\textbf{TS+V}($\kappa$)'' (from ``Tightening and Screening with Valid inequalities'').

To derive the set of valid inequalities for problem \eqref{MIPGenCC}, we define the real variables $z_j\in\mathbb{R}$ as $z_j:=\hat{a}_j^{\top} x$, and let $z_j^d=\inf_{x\in X}\hat{a}_j^{\top} x$,  $z_j^u=\sup_{x\in X}\hat{a}_j^{\top} x$ denote the lower and upper bounds on $z$ induced by the feasibility set \eqref{MIPGenCC_xinX}. Let us also define the function $L_{js} : f_{js}(z_j) = \Omega_{js} z_j - b_{js}$ for $z_j \in [z_j^d, z_j^u]$ and the set of functions $\mathcal{L}_j := \left\{ L_{js}, \; \forall s \in \mathcal{S} \right\}$. The valid inequalities we propose are heavily supported by the concepts of \emph{$k$-lower} and \emph{$k$-upper envelopes}, which we define in the following.

\begin{definition}
 For a given line $L_{js}$, we say that the point $(\tilde{z},\tilde{t})\in\mathbb{R}^2$ \emph{lies below}, \emph{on} or \emph{above} function $L_{js}$ depending on whether $\tilde{t}<\Omega_{js} \tilde{z} - b_{js}$, $\tilde{t}=\Omega_{js} \tilde{z} - b_{js}$ or $\tilde{t}>\Omega_{js} \tilde{z} - b_{js}$, respectively.
Naturally, we also say that function $L_{js}$ \emph{lies above}, \emph{contains}, or \emph{lies below} point $(\tilde{z},\tilde{t})$ in these cases. We also say that a point $(\tilde{z},\tilde{t})$ belongs to the set of lines $\mathcal{L}_j$ if there exists a line $L_{js}\in\mathcal{L}_j$ that contains the point $(\tilde{z},\tilde{t})$. 
\end{definition}
\begin{definition}
For a set of lines $\mathcal{L}_j$, the \emph{lower} (resp.\ \emph{upper}) \emph{score} of a point is the number of lines in $\mathcal{L}_j$ that lie below (resp.\ above) that point. The \emph{$k$-lower} (resp.\ \emph{$k$-upper}) \emph{envelope} of a set of lines $\mathcal{L}_j$ is the closure of the set of points that belong to $\mathcal{L}_j$ and that have lower (resp.\ upper) score equal to $k-1$. The $k$-lower envelope is also known as \emph{$k$-level}.
\end{definition}

\begin{figure}
\centering
\includegraphics[width=0.7\textwidth]{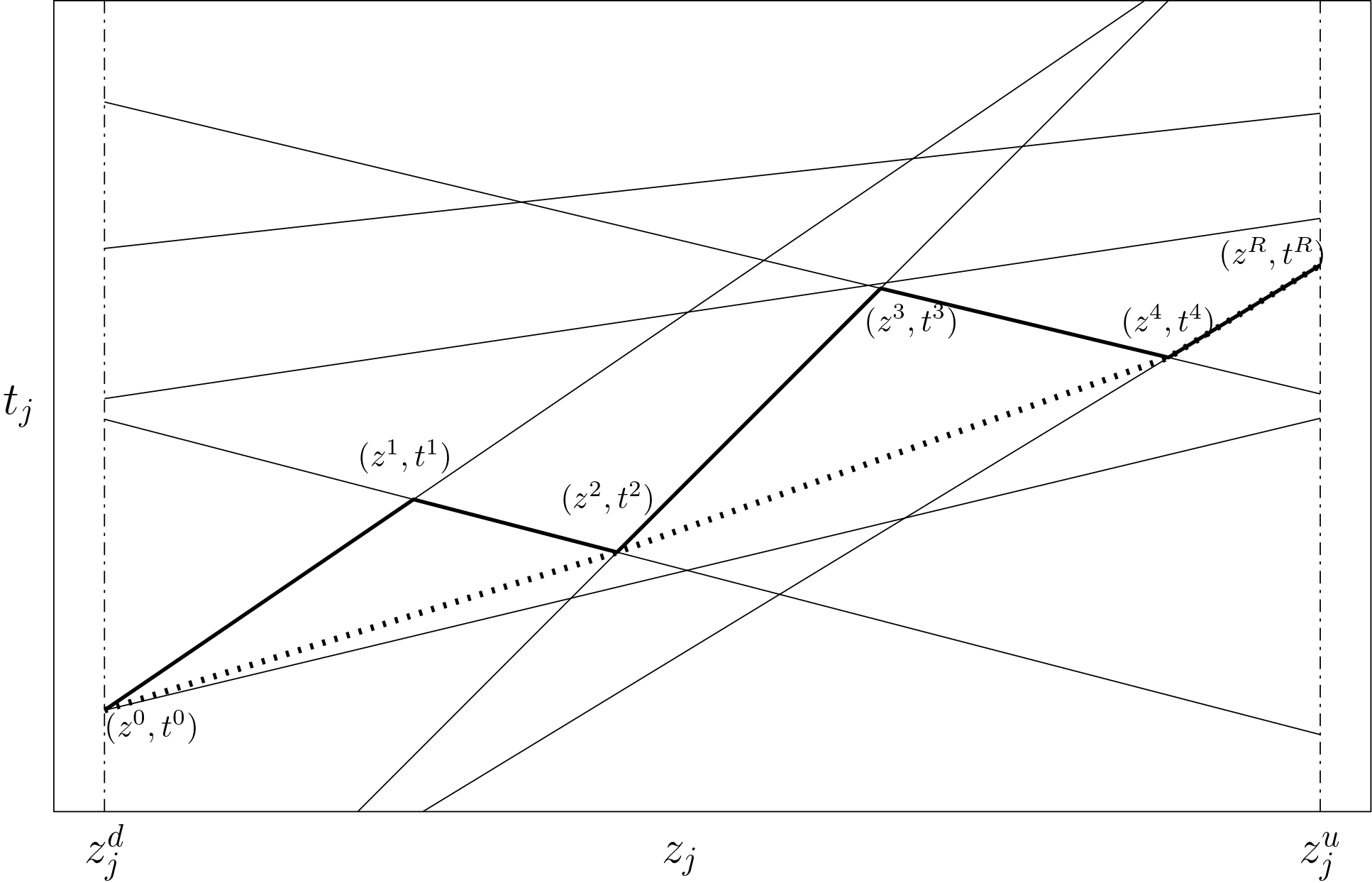} 
\caption{In bold, the 5-upper envelope (4-lower envelope, 4-level) of a set $\boldsymbol{\mathcal{L}_j}$ of 8 lines in the plane. In dotted, the lower hull of the 5-upper envelope}\label{fig:kenvelope} 
\end{figure}

For the sake of illustration, Figure \ref{fig:kenvelope} shows in bold the $5$-upper envelope of a set of 8 lines. Clearly, the $k$-envelopes of sets $\mathcal{L}_j$ can be seen as piece-wise linear functions on $z_j$.

\begin{proposition} \label{prop:VVII_versionpiecewise}
For a fixed $j\in \mathcal{J}$, let $U_j^{p+1}(\bigcdot)$ be the (p+1)-upper envelope of the set of lines $\mathcal{L}_j$, with $p:=  \floor*{\epsilon |\mathcal{S}|}$. Then the inequality
\begin{equation} \label{eq:VVII_versionpiecewise}
    U_j^{p+1}(\hat{a}_j^{\top} x) + x^{\top} a_j^{0} \le 0, \quad  x\in X 
\end{equation}
\noindent is valid for problem \eqref{MIPGenCC}.
\end{proposition}

\proof
Let $\bar{x},\bar{y}$ be any feasible solution of problem \eqref{MIPGenCC} with $\bar{z}_j = \hat{a}_j^{\top} \bar{x} \in [z_j^d,z_j^u]$, and suppose that $U^{p+1}_j(\bar{z}_j) + \bar{x}^{\top} a_j^{0} > 0$. By definition of  $k$-upper envelope, there exist $p+1$ lines in $\mathcal{L}_j$, $L_{js_i}: f_{js_i}(z_j)= \Omega_{js_i}z_j - b_{js_i}$ $\forall i\in \{1,\dots,p+1\}$ such that $f_{js_i}(\bar{z}_j)= \Omega_{js_i}\bar{z}_j - b_{js_i} \ge U_j^{p+1}(\bar{z}_j)$. 
Then for $i\in \{1,\dots,p+1\}$ it holds that
$\Omega_{js_i}\bar{z}_j - b_{js_i} + \bar{x}^{\top} a_j^{0} > 0$. Substituting in constraint \eqref{MIPGenCC_chance}, we obtain that $\bar{y}_{s_i}M_{s_i}>0$ $\Rightarrow$ $\bar{y}_{s_i}=1$, for $i\in \{1,\dots,p+1\}$. But then constraint \eqref{MIPGenCC_violation} is not satisfied, since $\sum_{s\in \mathcal{S}} \bar{y}_s \ge \sum_{i=1}^{p+1} \bar{y}_{s_i} = p+1 > p$. This is, however, in contradiction with our initial statement that $\bar{x}, \bar{y}$ is a feasible solution of problem \eqref{MIPGenCC}.
\endproof

The technical literature already includes references that propose methodologies to determine the $k$-envelope of a set of linear functions. For instance, the authors of \cite{cheema2014} give a basic algorithm for constructing $k$-envelopes called the Rider Algorithm. Essentially, this algorithm is based on the fact that the $k$-envelope is an unbounded polygonal chain that can be described by a sequence of vertices, which are intersections of lines of the set. In fact, every point in the $k$-envelope is contained in a given line. In this paper, we adapt the algorithm proposed in \cite{cheema2014} to the particular case in which $z_j^d,z_j^u$ are finite. Algorithm \ref{alg:Rider} describes our procedure in detail for a general set of the type $\mathcal{L}_j$.

\begin{algorithm}
\begin{small}
\caption{Rider Algorithm to build the $k$-upper envelope of $\mathcal{L}_j$ (adapted to the bounded case) } \label{alg:Rider}
\begin{algorithmic} 
\STATE To describe the $k$-upper envelope of a set of lines $\mathcal{L}_j$, we derive the sequence of vertices of the polygonal chain (intersections of lines in $\mathcal{L}_j$) $((z^0,t^0),\dots, (z^R, t^R))$ with $(z_j^d = z^0 < \dots<z^R=z_j^u)$. \\
\begin{enumerate}[\hspace*{0.5cm} Step 1.]\itemsep -.1cm
\item Set $r=0$. Let $\{(z^d_j,\Omega_{js}z_j^d-b_{js}),\forall s\in \mathcal{S}\}$, be the intersections of the lines $L_{js}\in\mathcal{L}_j$ with the vertical line $z=z_j^d$, and assume w.l.o.g.\ that all the points have different upper scores. Let $s_0\in \mathcal{S}$ be such that the intersection $(z^d_j,\Omega_{js_0}z_j^d-b_{js_0})$ has upper score equal to $k-1$. Then $(z^0,t^0)=(z^d_j,\Omega_{js_0}z_j^d-b_{js_0})$.
\item Compute the value of $z_j$ for which the line $L_{js_r}$ intersects the rest of the lines $L_{js}\in\mathcal{L}_j$ as $\text{int}_z(s,s_r) = \frac{b_{js}-b_{js_r}}{\Omega_{js}-\Omega_{js_r}}$. If $\exists$ $s'$ : $z^r < \text{int}_z(s',s_r) < z^u$, go to Step 3. Otherwise, go to Step~4.  
\item Find the line that intersects $L_{js_r}$ at the leftmost point to the right of $z^r$, that is, find $s_{r+1}=\arg\min_s \left\{\text{int}_z(s,s_r): \text{int}_z(s,s_r) > z^{r}\right\}$. Set $(z^{r+1},t^{r+1})=( \text{int}_z(s_{r+1},s_r),\Omega_{js_{r}}\text{int}_z(s_{r+1},s_r)-b_{js_{r}})$, update $r=r+1$ and go to Step 2. 
\item Set $(z^R,t^R)=(z^u_j,\Omega_{js_r}z_j^u-b_{js_r})$.
\end{enumerate}
\end{algorithmic}
\end{small}
\end{algorithm}
The LHSs of the valid inequalities \eqref{eq:VVII_versionpiecewise} are piece-wise linear functions not necessarily convex. Therefore, the inclusion of these inequalities into model \eqref{MIPGenCC} would require a significant amount of additional binary variables, which, in turn, is expected to increase the computational burden of this problem even further. Alternatively, we compute the lower hull of the $k$-upper envelope of $\mathcal{L}_j$, which takes the form of a convex piece-wise linear function. From this lower hull, we can extract a set of linear valid inequalities (the linear extensions of the pieces) that can be seamlessly inserted into model \eqref{MIPGenCC} without the need of any extra binary variables. In doing so, we are able to tighten model \eqref{MIPGenCC}, which can thus be solved more efficiently by available optimization software. Before presenting the set of valid inequalities, the following definition is required.

\begin{definition}
Let $Z$ be the convex hull of a set of points $P$. The \emph{upper} (resp.\ \emph{lower}) \emph{hull} of $P$ is the set of edges of $Z$ that lie on or above (resp.\ on or below) every point in $P$. 
\end{definition}

Corollary \ref{cor:VVIIconvexificadas} presents the set of linear valid inequalities \eqref{eq:VVII_versionupperhull} given by the lower hull of the $k$-upper envelope of~$\mathcal{L}_j$. 

\begin{corollary} \label{cor:VVIIconvexificadas}
Let $\{(z^r,t^r)\}$, $r\in \{0,\dots,R\}$, be the ordered set of vertices obtained by applying Algorithm \ref{alg:Rider} to set $\mathcal{L}_j$, and let $\{(z^{r'},t^{r'})\}$, $r'\in\{0,\dots,R'\} \subseteq \{0,\dots,R\}$, be the ordered subset of vertices such that the associated polygonal chain is the lower hull. Then the following linear inequalities 
\begin{equation} \label{eq:VVII_versionupperhull}
	\frac{t^{r'+1}-t^{r'}}{z^{r'+1}-z^{r'}}(\hat{a}_j^{\top} x - z^{r'}) + t^{r'} \le  - x^{\top} a_j^{0}, \quad x\in X, r'\in\{0,\dots,R'-1\}
\end{equation}
\noindent are valid for problem \eqref{MIPGenCC}.
\end{corollary}
\proof
The proof is straightforward, since for each $x\in X$ it holds $\frac{t^{r'+1}-t^{r'}}{z^{r'+1}-z^{r'}}(\hat{a}_j^{\top} x - z^{r'}) + t^{r'} \le U_j^{p+1}(\hat{a}_j^{\top} x) \le  - x^{\top} a_j^{0}$, by hypothesis and using \eqref{eq:VVII_versionpiecewise}.
\endproof

There exist plenty of algorithms of convexification of a set of points in the plane. Two of the most well-known are the \emph{Jarvis march} and the \emph{Graham scan} (\cite{toth2017}). Here, we give a simplified version of the former where we make use of special features of our set of points and only compute the lower hull. In particular, we assume that we have a set of presorted points $\{(z^{r},t^{r})\}$, $r\in\{0,\dots,R\}$ whose first and last points always belong to the hull. To speed up the process, we can find the point $(z^{\bar{r}}, t^{\bar{r}})$ from the previous set with the lowest $t$-coordinate (which always belongs to the lower hull), and then apply the algorithm to the subsets $\{(z^0,t^0),\dots,(z^{\bar{r}},t^{\bar{r}})\}$ and $\{(z^{\bar{r}},t^{\bar{r}}),\dots, (z^{R},t^{R})\}$. Algorithm  \ref{alg:JarvisMarch} describes in detail the proposed convexification procedure. 

\begin{algorithm}
\begin{small}
\caption{Jarvis March Algorithm to obtain the lower hull of a set of presorted points~$P$} \label{alg:JarvisMarch}
\begin{algorithmic} 
\STATE Let $P=\{(z^r,t^r)\}$, $r\in\{0,\dots,R\}$ be a set of points with $z^0 <\dots<z^{R}$. The algorithm derives a subset $P'=\{(z^{r'},t^{r'})\}$, $r'\in\{0,\dots,R'\}\subseteq \{0,\dots,R\}$, which constitutes the lower hull of the first set.\\
\begin{enumerate}[\hspace*{0.5cm} Step 1.]\itemsep -.1cm
\item Initially, set $P' = \{(z^0,t^0)\}$. 
\item Assume $(z^{r'},t^{r'})$ is the last point included in $P'$. If $r'=R$, the lower hull of $P$ is given by the set of points $P'$.
Otherwise, let $r'+1:= \arg \min_r\left\{	\frac{t^{r}-t^{r'}}{z^{r}-z^{r'}}: (z^r,t^r)\in P \text{ with } z^r>z^{r'} \right\}$. Update $P' = P' \cup \{(z^{r'+1},t^{r'+1})\}$. Repeat Step 2.
\end{enumerate}
\end{algorithmic}
\end{small}
\end{algorithm}

Algorithm \ref{alg:Rider} constructs the $k$-upper envelope of $\mathcal{L}$ in $O(n_k\log^2 |\mathcal{L}|)$ time and $O(n_k|\mathcal{L}|+|\mathcal{L}|)$ space, where $n_k$ denotes the maximum number of edges of the $k$-upper envelope of any set with $|\mathcal{L}|$ lines (see \cite{edelsbrunner1986}). An upper bound on $n_k = O(|\mathcal{L}|k^{\frac{1}{2}})$ can be found in \cite{edelsbrunner1985}. \cite{dey1998,toth2000,nivasch2008} have subsequently improved the upper bound of the complexity of this problem by studying a closely related problem, the $k$-set problem, in the two-dimensional case. As for the Jarvis March Algorithm, in the two-dimensional case presented here it has running time $O((R+2)(R'+1))$ (\cite{toth2017}). For detailed results on the complexity of related algorithms and for extensions to the multidimensional case, we refer the reader to  \ref{sec:anexoExtensionMultidim}. 

Valid inequalities \eqref{eq:VVII_versionpiecewise} and Algorithms \ref{alg:Rider} and \ref{alg:JarvisMarch} can be generalized to work with finite discrete distributions with unequal probabilities. Observe that the definition of the $(p+1)$-lower (resp.\ $(p+1)$-upper) envelope is based on satisfying constraint \eqref{MIPGenCC_violation}. In this extension, constraint \eqref{MIPGenCC_violation} is replaced with $\sum_{s\in \mathcal{S}} p_sy_s \le \epsilon$, where $p_s$ is the probability of scenario $s$. Hence, it suffices to redefine the lower (resp.\ upper) score of a point as the \textit{probability of the} lines that lie below (resp.\ above) that point, and modify the definition of $k$-lower (resp.\ upper) envelope and Step 3\ of Algorithm~\ref{alg:Rider} accordingly, so that this more general constraint is satisfied.

The following theorem sheds light on the strength of our valid inequalities.
\begin{theorem}
Let $\mathcal{FR}_j$ be the feasible region of the relaxation of problem \eqref{MIPGenCC} associated to a given row $j\in \mathcal{J}$, i.e.,
\begin{subequations}
\label{MIPGenCCj}
\begin{align}
 \min_{x,y_s} \quad &  f(x) \label{MIPGenCCj_FO}\\
\text{s.t.}  \quad & \eqref{MIPGenCC_xinX}, \eqref{MIPGenCC_violation}, \eqref{MIPGenCC_ybinaria}  \label{MIPGenCCj_resto}\\
& \Omega_{js} \hat{a}_j^{\top} x - b_{js} + x^{\top} a_j^{0} \le M_{js}y_s, \quad \forall s\in \mathcal{S}. \label{MIPGenCCj_chance}
\end{align}
\end{subequations}

Assuming that the set $X$ of deterministic constraints given in \eqref{MIPGenCC_xinX} is non-empty and compact, the set \eqref{eq:VVII_versionupperhull} of valid inequalities for all $j \in \mathcal{J}$ provides $\bigcap_{j\in \mathcal{J}} \mathrm{conv}(\mathcal{FR}_j)$.
\end{theorem}

The proof of this theorem is a direct consequence of Proposition~\ref{prop:UequalFR} in \ref{sec:anexoQuantileCuts}.

Furthermore, it is possible to relate our valid inequalities with the \emph{quantile cuts} described in \cite{qiu2014} and \cite{xie2018}. To see this, consider problem \eqref{MIPGenCC} with feasible region $\mathcal{FR}$ and note that $\mathcal{FR} = \bigcap_{j\in \mathcal{J}}\mathcal{FR}_j$. The \textit{quantile cuts} can be viewed as a projection of the well-known family of mixing inequalities in the $(x, y)$-space onto the original space $x$. Quantile cuts represent an infinite family of inequalities, and the intersection of all quantile cuts is called the \textit{quantile closure}. If another round of quantile cuts is applied to a stronger formulation, specifically the original formulation with the (first) quantile closure, a subsequent closure can be derived (the second closure), and so on. In \cite{xie2018} it is proved that the sequence of sets obtained by successive quantile closure operations converges to the convex hull of the feasible region of the original problem with respect to the Hausdorff metric, $\mathrm{conv}(\mathcal{FR})$. It is also shown that the separation over the first quantile closure is NP-hard. Corollary~\ref{cor:clausuramedianteVVII} in \ref{sec:anexoQuantileCuts} shows that the set \eqref{eq:VVII_versionupperhull} of valid inequalities (per constraint $j \in \mathcal{J}$) coincides with the limit of the succession of quantile closures of $\mathcal{FR}_j$ as the order of the closure grows to infinity. This limit is precisely $\mathrm{conv}(\mathcal{FR}_j)$.

\section{Optimal Power Flow under Uncertainty: A Joint Chance-Constrained Modeling Approach} \label{sec:OPF}
The OPF is a routine at the core of important tools for power system operations (\cite{frank2012optimal}). In its deterministic version, the OPF problem seeks to determine the least-costly dispatch of thermal generating units to satisfy the system's \emph{net} demand (i.e., demand minus renewable generation), while complying with the technical limits of production and transmission network equipment. The main challenge of the OPF is that it is a non-linear and non-convex optimization problem, due to the power flow equations that govern the (static) behavior of power systems. For this reason, the \emph{direct current} approximation (DC) of the power flow equations, which transforms the problem into a linear program, is frequently used. The demand and renewable generation are factors that increase the uncertainty in power systems, and ignoring it can lead to unsafe operating conditions. 

 However, given the inherently uncertain nature of the electricity net demand, the probabilistic version of the OPF problem can be formulated as a JCCP that aims at minimizing the expected production cost while enforcing that the technical constraints are satisfied with a given (high) probability. In this context, chance-constrained programming can be used to minimize the expected operating cost whilst guaranteeing that the system withstands unforeseen peeks of electrical load due to stochastic demand or uncertainty in power generation (\cite{vanackooij2011}). The chance-constrained OPF (CC-OPF) problem addresses this uncertainty and pursues to ensure the safe operation of a power system with a high level of probability.

To address the CC-OPF problem, several papers in the literature (e.g., \cite{LineGoran}) directly work with SCCs. However, the main drawback of this modeling approach is that, even in those cases where the probability of violating each individual constraint seems more than tolerable, the resulting \emph{joint risk} (that is, the probability that \emph{any} of the technical constraints be violated) may still be excessive and inadmissible. This is the key motivation behind the use of JCCs to tackle the CC-OPF problem (see, e.g.,  \cite{LineAlejandra}).

In this vein, there are several approaches in the literature to solve the JCC-OPF problem. \cite{vrakopoulou2013} adopt the scenario approach (SA) to approximate the solution of the JCC-OPF, while \cite{chen2021time} propose a heuristic data-driven method that involves enforcing the satisfaction of the technical constraints for a box of the uncertainty. This box is inferred using one-class support vector clustering and its size is contingent on the system's desired reliability. \cite{LineTunning} propose an iterative tuning algorithm to solve a robust reformulation of the JCC-OPF problem. \cite{esteban2023distributionally} introduce a distributionally robust JCC-OPF model that considers contextual information using an ambiguity set based on probability trimmings. To make their model tractable, they resort to the widely known CVaR-based approximation of the JCC. 

The aforementioned SAA method is another effective way to solve JCCPs and has the potential to identify OPF solutions with a better cost performance than that of the more conservative solutions delivered by the previous approaches. However, solving the JCC-OPF problem using SAA is challenging due to the presence of binary variables, the number of scenarios required and the size of the power systems. \cite{lejeune2020optimal} propose a methodology to solve the SAA of the JCC-OPF without including the power flow equations into the joint chance-constrained system. To the best of our knowledge, we are the first to efficiently solve the JCC-OPF problem by means of the SAA approach, using a MIP reformulation and including the arduous power flow constraints. Furthermore, unlike the sample-based approach introduced in \cite{LineAlejandra}, which relies on a smooth nonlinear approximation of the JCC-OPF, ours offers optimality guarantees.

\subsection{Notation and problem formulation} \label{sec:OPF_notation}
In this subsection, we present the  formulation of the JCC-OPF problem that we consider in the case study. To do so, we introduce the following notation and modeling choices, widely used in the power system domain:
\begin{enumerate}
    \item \emph{Power system:} A power system can be represented by a directed graph denoted by $G(\mathcal{N},\mathcal{L})$, where $\mathcal{N}$ is the set of buses (vertices), indexed by $n$, and $\mathcal{L}$ is the set of transmissions lines (edges), whose direction is chosen arbitrarily and are indexed by $l$. To supply the net demand at each bus, the power system is equipped with a set of generators $\mathcal{G}$, indexed by $g$, where $\mathcal{G}_n$ represents the set of generators at bus $n$. Note that the energy injected at each bus, derived from the balance between generation and net demand, is transported among the buses through the transmission lines. For illustration, Figure \ref{fig:3bus} depicts a representation of a power system with three buses, three transmission lines, two generators located at nodes 1 and 2 and one electricity net demand at node 3.
    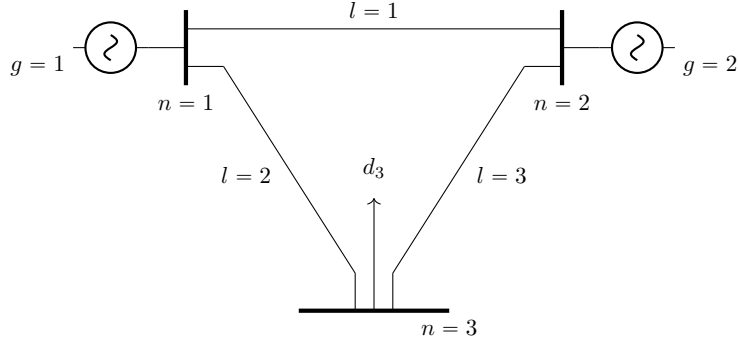
\begin{figure}
    \centering
    \begin{circuitikz}[scale=1,every node/.style={scale=0.8}]
    \draw [ultra thick] (0,4)  -- (0,3) node[anchor=north]{$n = 1$};
    \draw [ultra thick] (5,4)  -- (5,3) node[anchor=north]{$n = 2$};
    \draw [ultra thick] (1.5,0)  -- (3.5,0) node[anchor=north]{$n= 3$};
    \draw[] (-1.5,3.5) node[anchor=north east]{$g = 1$} to [sV] (-0.5,3.5);
    \draw (-0.5,3.5) -- (0,3.5);
    \draw[] (6.5,3.5) node[anchor=north west]{$g = 2$} to [sV] (5.5,3.5); 
    \draw (5.5,3.5) -- (5,3.5); 
    \draw (0,3.75) -- (5,3.75);
    \node[] at (2.5,4) {$l = 1$ };
    \draw (0,3.25) -- (0.5,3.25);
    \draw (0.5,3.25) -- (2.25,0.5);
    \draw (2.25,0.5) -- (2.25,0);
    \node[] at (0.8,1.8) {$l = 2$};
    \draw (2.75,0) -- (2.75,0.5);
    \draw (2.75,0.5) -- (4.5,3.25);
    \draw (4.5,3.25) -- (5,3.25);
    \node[] at (4.2,1.8) {$l = 3$};
    \draw[->] (2.5,0.0) -- (2.5,1.5); 
    \node[] at (2.5,1.9) {$d_3$};
    \end{circuitikz}
    \caption{Illustrative three-node power system}
    \label{fig:3bus}
    \end{figure}
    
    \item \emph{Nodal net loads}: The (uncertain) electricity net demand at node $n$, $\tilde{d}_n$, is given by $\tilde{d}_n=d_n - \omega_n$, where $d_n$ is the predicted value and $\omega_n$ is the  forecast error with a change of sign. This error is modeled as a random variable with zero mean which follows an unknown continuous probability distribution.
    \item \emph{Generation}: To cope with the forecast errors $(\omega_n)_{n \in \mathcal{N}}$, generators' power outputs are adjusted according to the following affine control policy:
    \begin{align*}
        & \tilde{p}_g = p_g - \beta_g\Omega, \quad \forall g \in \mathcal{G},
    \end{align*}
    where $\Omega:= \sum_{n \in \mathcal{N}} \omega_n$ is the system-wise aggregated forecast error, and $p_g$ and $\beta_g$ are the power output dispatch and the participation factor of generating unit $g$, respectively (see, e.g., \cite{bienstock2014chance,LineTunning,LineAlejandra}). The minimum and maximum capacity of generator $g$ is denoted by $\underline{p}_g$ and $\overline{p}_g$, respectively.
    \item \emph{Power balance}: Given the affine control policy of the previous point, the power balance equation takes the following form:
    \begin{align*}
        & \sum_{g \in \mathcal{G}} \tilde{p}_g  - \sum_{n \in \mathcal{N}} \tilde{d}_n  = \sum_{g \in \mathcal{G}} \left(p_g - \beta_g\Omega\right) - \sum_{n \in \mathcal{N}} \left(d_n - \omega_n\right) = 0. 
    \end{align*}
    Hence, to ensure the power balance for \emph{any} realization of the forecast errors $(\omega_n)_{n \in \mathcal{N}}$, it must hold:
    \begin{eqnarray*}
    \sum_{g \in \mathcal{G}} p_g - \sum_{n \in \mathcal{N}} d_n &=& 0\\
    \sum_{g \in \mathcal{G}} \beta_g &=& 1.
    \end{eqnarray*}
    \item \emph{Power flows}: Line flows are modeled using the well-known approximation based on the power transfer distribution factors (PTDFs), $B_{ln}$, $l \in \mathcal{L}$, $n \in \mathcal{N}$,   which sets a linear relation between the power flow through line $l$ and the power injected at node $n$. The maximum capacity of line $l$ is denoted by $\overline{f}_l$.
    \item \emph{Power production cost}: The cost function of each generating unit is assumed to be quadratic and, as a result, the total power production cost is given by
    \begin{align*}
        & \sum_{g \in \mathcal{G}} C_{2,g} \, \left(p_{g} - \Omega\beta_g\right)^2 + C_{1,g} \, \left(p_{g} - \Omega\beta_g\right) + C_{0,g},
    \end{align*}
    \noindent where $C_{2,g}$, $C_{1,g}$, $C_{0,g}$ are the  coefficients defining the quadratic cost function of generating unit $g$. On the assumption that $\omega_n$, for each $n \in \mathcal{N}$, is a random variable with zero mean, we have (see, for instance, \cite{LineAlejandra})
    \begin{align*}
        & \mathbb{E} \left[ \sum_{g \in \mathcal{G}} C_{2,g} \, \left(p_{g} - \Omega\beta_g\right)^2 + C_{1,g} \, \left(p_{g} - \Omega\beta_g\right) + C_{0,g} \right] = \sum_{g \in \mathcal{G}} C_{2,g} \, p_{g}^2 + C_{1,g} \, p_{g}+ C_{0,g} + \mathbb{V} \left(\Omega\right) \, C_{2,g} \, \beta_{g}^2,
    \end{align*}
    \noindent where $\mathbb{V}(\Omega)$ denotes the variance of the random variable $\Omega$.
    
\end{enumerate}

With the above ingredients, the JCC-OPF problem that we tackle in this paper is formulated as follows:
\begin{subequations}
\label{eq:JCC-OPF}
\begin{align}
\min_{p_g,\beta_g \geq 0, \forall g \in \mathcal{G}} \quad & \sum_{g \in \mathcal{G}} C_{2,g} \, p_{g}^2 + C_{1,g} \, p_{g}+ C_{0,g} + \mathbb{V} \left(\Omega\right) \, C_{2,g} \, \beta_{g}^2 \label{eq:OPF_objective}\\
\text{s.t.} \quad & \sum_{g \in \mathcal{G}} \beta_g = 1 \label{eq:OPF_balance1}\\
& \sum_{g \in \mathcal{G}} p_{g} -  \sum_{n \in \mathcal{N}} d_{n} = 0 \label{eq:OPF_balance2}\\
&\underline{p}_{g} \leq p_g \leq \overline{p}_{g}, \quad \forall g \in \mathcal{G} \label{eq:OPF_gen-det}\\
& -\overline{f}_{l} \leq \sum_{n \in \mathcal{N}} B_{ln}\left(\sum_{g \in \mathcal{G}_n} p_{g} - d_n \right) \leq \overline{f}_{l}, \quad \forall l \in \mathcal{L} \label{eq:OPF_flow-det}\\
& \mathbb{P}
\left(\begin{array}{l}
    \underline{p}_{g} \leq p_g -\Omega\beta_g \leq \overline{p}_{g}, \quad \forall g \in \mathcal{G} \\
    -\overline{f}_{l} \leq \displaystyle\sum_{n \in \mathcal{N}} B_{ln}\left(\displaystyle\sum_{g \in \mathcal{G}_n} \left(p_{g} - \Omega\beta_g\right) + \omega_n - d_n \right) \leq \overline{f}_{l}, \quad \forall l \in \mathcal{L} 
\end{array} \right) \geq 1 - \epsilon \label{eq:OPF_jointCC}.
\end{align}
\end{subequations}

The objective \eqref{eq:OPF_objective} is the minimization of the expected total generation cost. The equality constraints \eqref{eq:OPF_balance1} and \eqref{eq:OPF_balance2} enforce the power balance in the system, whereas constraints \eqref{eq:OPF_gen-det} and \eqref{eq:OPF_flow-det} ensure a feasible power dispatch which  corresponds to an error-free scenario, i.e., to a realization of the net-load forecast errors such that $\omega_n = 0$ $\forall n \in \mathcal{N}$. Finally, expression \eqref{eq:OPF_jointCC} constitutes the joint chance-constraint system by which the decision-maker states that the OPF solution must be feasible with a probability greater than or equal to $1-\epsilon$. Accordingly, parameter $\epsilon$ is the maximum allowed probability of constraint violation set by the user. Formulation \eqref{eq:JCC-OPF} is quite standard and has been used before by \cite{bienstock2014chance,LineTunning,LineAlejandra}, among others.

Problem~\eqref{eq:JCC-OPF} can be written in the form of~\eqref{GenCC}, i.e., as a CCP with linear JCC and random RHS and LHS. To see this, define the vector of continuous decision variables $x$ in \eqref{GenCC} as $x:=(p_g,\beta_g)_{g \in \mathcal{G}}$, and group all these variables by means of the set $\mathcal{I}$ with elements $i$ running from 1 to $|\mathcal{I}| = 2|\mathcal{G}|$. In this way, we have that  $x \in \mathbb{R}^{|\mathcal{I}|}_{+}$, the set $X\subseteq \mathbb{R}^{|\mathcal{I}|}$ represents the polyhedron defined by the deterministic constraints \eqref{eq:OPF_balance1}--\eqref{eq:OPF_flow-det}, and $f: \mathbb{R}^{|\mathcal{I}|} \longrightarrow \mathbb{R}$ is  the convex function providing the expected total generation cost \eqref{eq:OPF_objective}. Likewise, if we collect all the constraints involved in the joint chance-constraint system \eqref{eq:OPF_jointCC} into the set $\mathcal{J}$ (hence $|\mathcal{J}| = 2|\mathcal{G}| + 2|\mathcal{L}|$), this system can be represented by way of a technology matrix with random rows $a_{j}(\omega)\in \mathbb{R}^{|\mathcal{I}|}$, $j \in \mathcal{J}$, and random RHS  $b_j(\omega) \in \mathbb{R}$, $j \in \mathcal{J}$, where the uncertainty is again represented through the random vector $\omega$ taking values in $\mathbb{R}^{|\mathcal{N}|}$. 

As discussed in Section \ref{sec:JCCP}, the CCP \eqref{eq:JCC-OPF} can be easily reformulated into a MIP using SAA. Thus, we assume that  $\omega$ has a finite discrete support defined by a collection of points $\{\omega_s \in \mathbb{R}^{|\mathcal{N}|}, s \in \mathcal{S}\}$ and respective probability masses $\mathbb{P}(\omega = \omega_s)=\frac{1}{|\mathcal{S}|}$, $\forall s\in \mathcal{S}=\{1,\dots,|\mathcal{S}|\}$. Accordingly, $\omega_{ns}$ and $\Omega_s$ are realizations of the respective random variables under scenario $s$. We define $p:=  \floor*{\epsilon |\mathcal{S}|}$, the vector $y$ of binary variables $y_s$, $\forall s\in \mathcal{S}$, and the large enough constants $M^1_{gs}, M^2_{gs}, M^3_{ls}, M^4_{ls}$. Thus, the MIP reformulation of problem \eqref{eq:JCC-OPF} writes as follows: 
\begin{subequations}
\label{eq:BigM-OPF}
\begin{align}
\min_{p_g,\beta_g\geq0,y_s} \quad & \sum_{g \in \mathcal{G}} C_{2,g} \, p_{g}^2 + C_{1,g} \, p_{g}+ C_{0,g} + \widehat{\mathbb{V}} \left(\Omega\right) \, C_{2,g} \, \beta_{g}^2  \label{eq:MIP1_objective}\\
\text{s.t.}  \quad & \eqref{eq:OPF_balance1} - \eqref{eq:OPF_flow-det} \label{eq:MIP1_repeated}\\
& p_g - \Omega_s \beta_g  \geq \underline{p}_{g} - y_s M^1_{gs}, \quad \forall g \in \mathcal{G}, s \in \mathcal{S} \label{eq:MIP1_gen_LB}\\
& p_g -\Omega_s \beta_g \leq \overline{p}_{g} + y_s M^{2}_{gs}, \quad \forall g \in \mathcal{G}, s \in \mathcal{S} \label{eq:MIP1_gen_UB}\\
& \sum_{n \in \mathcal{N}} B_{ln}\left(\sum_{g \in \mathcal{G}_n} \left(p_{g} - \Omega_s \beta_g\right) - d_n + \omega_{ns} \right) \geq -\overline{f}_{l} - y_s M^{3}_{ls}, \quad \forall l \in \mathcal{L}, s \in \mathcal{S} \label{eq:MIP1_flow_LB}\\
& \sum_{n \in \mathcal{N}} B_{ln}\left(\sum_{g \in \mathcal{G}_n} \left(p_{g} - \Omega_s \beta_g\right) - d_n + \omega_{ns} \right) \leq \overline{f}_{l} + y_s M^4_{ls}, \quad \forall l \in \mathcal{L}, s \in \mathcal{S} \label{eq:MIP1_flow_UB}\\
& \sum_{s \in \mathcal{S}} y_s \leq p \label{eq:MIP1_violation}\\
& y_s \in \{0,1\}, \quad \forall s \in \mathcal{S}. \label{eq:MIP1_bincharacter}
\end{align}
\end{subequations}

Constraints \eqref{eq:MIP1_gen_LB}-\eqref{eq:MIP1_flow_UB} represent the sample-based reformulation of the joint chance constraint \eqref{eq:OPF_jointCC}. For a given scenario $s\in \mathcal{S}$, inequalities \eqref{eq:MIP1_gen_LB}-\eqref{eq:MIP1_flow_UB} guarantee that all the original constraints are satisfied when $y_s = 0$. If $y_s=1$, some of the original constraints can be violated for scenario $s$. Finally, the inequality \eqref{eq:MIP1_violation} ensures that the probability of the JCC is met and the binary character of variables $y_s$ is declared in \eqref{eq:MIP1_bincharacter}. 

\subsection{Numerical Experiments} \label{sec:OPF_casestudy}
This section discusses a series of numerical experiments with which we evaluate the different approaches presented in Section \ref{sec:JCCP} to solve the SAA-based MIP reformulation of the JCC-OPF. In particular, we compare the performance of approaches \textbf{T}, \textbf{TS} and \textbf{TS+V} using five standard power systems widely employed in the technical literature on the topic, namely the IEEE-RTS-24, IEEE-57, IEEE-RTS-73, IEEE-118, and IEEE-300 test systems. All data pertaining to these systems are publicly available in the repository \cite{pglib} under version~21 and their main features are listed in Table~\ref{tab:features}. All optimization problems have been solved using GUROBI 9.1.2 (\cite{gurobi}) on a Linux-based server with CPUs clocking at 2.6 GHz, 6 threads and 32 GB of RAM. In all cases, the optimality GAP has been set to $10^{-9}\%$ and the time limit to 10 hours. All data used and codes implemented throughout this numerical experiments can be found in the repository \cite{TCSAAJCCOPF2023}.

\begin{table}[htbp]
\begin{center}
\caption{Short Description of Test Power Systems}
\begin{tabular}{lccccc}
\hline
 &IEEE-RTS-24  &IEEE-57 &IEEE-RTS-73 &IEEE-118 &IEEE-300\\
\hline
\# Nodes &24 &57 &73 &118 &300\\
\# Generators &32 &4 &96 &19 & 57\\
\# Lines &38 &41 &120 &186 &411\\
\hline
\end{tabular}
\label{tab:features}
\end{center}
\end{table}

Similarly to \cite{LineAlejandra}, we assume that the error of net loads is normally distributed, i.e., $\omega \sim N(\mathbf{0},\Sigma)$, where $\mathbf{0}$ and $\Sigma$ represent, respectively, the zero vector and the covariance matrix. We also assume that the standard deviation of $\omega_n$ at node $n$ is proportional to the net nodal demand $d_n$ according to a parameter $\zeta$ between 0 and 1. Thus, this parameter controls the magnitude of net demand fluctuations. Under these assumptions, the procedure to model uncertainty proposed in \cite{LineAlejandra} runs as follows. First, we compute the positive definite matrix $C=\widehat{C}\widehat{C}^{\top}$ where each element of matrix $\widehat{C}$ is a sample randomly drawn from a uniform distribution with support in $[-1,1]$. Then, to obtain a positive definite matrix $\Sigma$ in which the diagonal elements are equal to $(\zeta d_{n})^{2}$, we define each of its entries ($\sigma_{nn'}$) as follows:
\begin{align*}
    & \sigma_{nn'} = \zeta^{2} \frac{c_{nn'}}{\sqrt{c_{nn}c_{n'n'}}} d_n d_{n'}, \quad \forall n,n' \in \mathcal{N},
\end{align*}
\noindent where $c_{nn'}$ denotes the element of matrix $C$ in row $n$ and column $n'$. To avoid generating infeasible instances of the JCC-OPF problem, the parameter $\zeta$ has been set to 0.15 for the four smallest systems and to 0.05 for the IEEE-300 system. To characterize the net demand uncertainty, we consider 1000 scenarios and a tolerable probability of violation of the joint chance constraint of 5\% (i.e., $\epsilon = 0.05$ and $p=50$). Finally, each solution strategy is run for ten different sets of randomly generated samples. Accordingly, in this section we provide tables with figures averaged over these ten instances. Also, the generated samples can be downloaded from the repository \cite{TCSAAJCCOPF2023}.

Table \ref{tab:bench} includes the results of solving the mixed-integer quadratic optimization model \eqref{eq:BigM-OPF} if the large constants $M^1_{gs}$, $M^2_{gs}$, $M^3_{ls}$ and $M^4_{ls}$ are set to a high enough arbitrary value, specifically  $10^4$. Despite being remarkably computationally expensive, this approach has been used in the technical literature (e.g., \cite{zhang2015data}), and thus we refer to it as \emph{benchmark approach} (\textbf{BN}). Table~\ref{tab:bench} includes the number of constraints in the model (\#CON), the linear relaxation gap (LR-GAP) calculated using the optimal solution of each instance, the optimality gap given by the difference between the best lower bound and the best integer solution found by the MIP solver (MIP-GAP), the number of instances solved to global optimality in less than 10 hours (\#OPT) and the solution time in seconds (Time). As expected, the computational time needed to solve the OPF with the Big-M model \eqref{eq:BigM-OPF} increases significantly with the size of the instances. While the 10 instances from systems IEEE-RTS-24, IEEE-57 and IEEE-73 are solved in less than 10 hours, none of the instances for systems IEEE-118 and IEEE-300 are solved to global optimality within that time limit, and the average MIP-GAP after 10 hours amounts to 0.29\% and 0.27\%, respectively. Interestingly, despite the fact that the LR-GAP is relatively low for the IEEE-RTS-73 system, the average computational time required to solve this case is particularly high compared to the two smaller systems.

\begin{table}[htbp]
\begin{center}
\caption{Benchmark approach (\textbf{BN}): Results}
\begin{tabular}{cccccc}
\hline
 &IEEE-RTS-24  &IEEE-57 &IEEE-RTS-73 &IEEE-118 &IEEE-300\\
\hline
\#CON  &140143	&168171	&432435	&410413 &936939\\
LR-GAP  &1.756\%	&0.623\%	&0.061\%	&0.956\%	&1.114\%\\
MIP-GAP (\#OPT) &0.00\% (10) &0.00\% (10) &0.00\% (10) &0.29\% (0) &0.27\% (0)\\
Time (s)  &1121.3	&103.2	&11161.2	&36000.0	&36000.0\\
\hline
\end{tabular}
\label{tab:bench}
\end{center}
\end{table}

As discussed in the technical literature, a proper tuning of the Big-Ms makes model \eqref{eq:BigM-OPF} tighter and generally easier to solve (\cite{qiu2014}) by the MIP routine. Therefore, we evaluate the computational performance of the ``Iterative Coefficient Strengthening'' Algorithm  and provide the corresponding results in Table \ref{tab:tightening}. In particular, \textbf{T}(1), \textbf{T}(2) and \textbf{T}(3) represent the results obtained by solving model \eqref{eq:BigM-OPF} with the Big-M values provided by Algorithm~\ref{alg:IterativeCoefStreng} with $\kappa=1$, 2 and 3, respectively. 
Table \ref{tab:tightening} includes the average values of LR-GAP and MIP-GAP, the number of instances solved to optimality in less than 10 hours (\#OPT) and the speedup factor with respect to the benchmark approach. To determine this factor, we have considered that the total computational time of approach \textbf{T} is given as the sum of the time required to run Algorithm \ref{alg:IterativeCoefStreng} $\kappa$ times to determine the Big-Ms plus the time it takes to solve problem \eqref{eq:BigM-OPF}. 

\begin{table}[htbp]
\begin{center}
\caption{Coefficient tightening approach (\textbf{T}): Results}
\begin{tabular}{ccccccc}
\hline
\multicolumn{2}{c}{} &IEEE-RTS-24  &IEEE-57 &IEEE-RTS-73 &IEEE-118 &IEEE-300\\
\hline
\multirow{3}{*}{\small LR-GAP} &\textbf{T}(1) &1.755\%	&0.510\%	&0.061\%	&0.711\%	&0.472\%\\
&\textbf{T}(2) &1.662\%	&0.330\%	&0.055\%	&0.522\%	&0.324\%\\
&\textbf{T}(3) &1.386\%	&0.255\%	&0.029\%	&0.434\%	&0.264\%\\
\hline
\multirow{3}{*}{\small MIP-GAP (\#OPT)} &\textbf{T}(1) &0.00\% (10) &0.00\% (10) &0.00\% (10) &0.27\% (0) &0.16\% (0)\\
&\textbf{T}(2) &0.00\% (10) &0.00\% (10) &0.03\% (2) &0.16\% (0) &0.09\% (0)\\
&\textbf{T}(3) &0.00\% (10) &0.00\% (10) &0.00\% (10) &0.12\% (0) &0.07\% (0)\\
\hline
\multirow{3}{*}{Speedup factor} &\textbf{T}(1) &0.22x	&0.07x	&0.61x &1.00x	&1.00x\\
&\textbf{T}(2) &0.17x	&0.13x	&0.31x	&1.00x	&1.00x\\
&\textbf{T}(3) &0.48x	&0.24x	&0.66x	&1.00x	&1.00x\\
\hline
\end{tabular}
\label{tab:tightening}
\end{center}
\end{table}

Since reducing the Big-M values makes model \eqref{eq:BigM-OPF} tighter, the results in Table \ref{tab:tightening} show lower values of LR-GAP with respect to  \textbf{BN}. Furthermore, this effect grows with the number of iterations since Algorithm~\ref{alg:IterativeCoefStreng} ensures that the Big-Ms never increase between iterations. Although decreasing the values of the Big-Ms leads to tighter MIPs for all the test systems, the numerical results in Table~\ref{tab:tightening} clearly indicate that computational savings are not guaranteed in all cases. Indeed, while the ten instances are solved by  \textbf{BN} in less than 10 hours for the IEEE-RTS-73 system, \textbf{T}(2) only provides the optimal solution for two instances. On top of that, the speedup factors for the three smaller systems are always lower than 1, which means that the computational times actually increase in these cases. On the contrary, the average MIP-GAP of the two largest systems is significantly decreased with respect to \textbf{BN}. Therefore, we conclude that, due to the heuristics implemented in current commercial MIP solvers, the computational advantages that one could expect a priori from ``Iterative Coefficient Strengthening'' are not always guaranteed and are contingent on the structure and data of the problem. 

Next, in Table~\ref{tab:screening} we provide the computational results related to the \textbf{TS} method, in which Algorithm~\ref{alg:IterativeCoefStreng} is extended to remove redundant constraints from model \eqref{eq:BigM-OPF}. Here, \#CON is provided as the percentage of the number of constraints of the reference model \textbf{BN} (indicated in Table~\ref{tab:bench}) that are retained by \textbf{TS} in each iteration. Table \ref{tab:screening} also includes the average MIP-GAP, the number of instances solved to optimality and the average speedup factor in relation to \textbf{BN}. 

\begin{table}[htbp]
\begin{center}
\caption{Tightening and Screening (\textbf{TS}): Results}
\begin{tabular}{ccccccc}
\hline
\multicolumn{2}{c}{} &IEEE-RTS-24  &IEEE-57 &IEEE-RTS-73 &IEEE-118 &IEEE-300\\
\hline
\multirow{3}{*}{\#CON}  &\textbf{TS}(1) &23.9\%	&2.8\%	&26.0\%	&8.9\%	&12.7\%\\
&\textbf{TS}(2) &23.3\%	&2.2\% &23.8\%	&6.3\%	&9.1\%\\
&\textbf{TS}(3) &23.1\%	&2.1\% &23.0\%	&5.6\%	&8.0\%\\
\hline
\multirow{3}{*}{\small MIP-GAP (\#OPT)} &\textbf{TS}(1) &0.00\% (10) &0.00\% (10) &0.00\% (10) &0.15\% (0) &0.08\% (0)\\
&\textbf{TS}(2) &0.00\% (10) &0.00\% (10) &0.00\% (10) &0.03\% (2) &0.04\% (0)\\
&\textbf{TS}(3) &0.00\% (10) &0.00\% (10) &0.00\% (10) &0.01\% (6) &0.01\% (4)\\
\hline
\multirow{3}{*}{Speedup factor} &\textbf{TS}(1) &1.5x	&1.8x	&4.7x	&1.0x	&1.0x\\
&\textbf{TS}(2) &1.5x	&3.8x	&2.6x	&1.1x	&1.0x\\
&\textbf{TS}(3) &3.8x	&4.4x	&15.1x	&1.4x	&1.2x\\
\hline
\end{tabular}
\label{tab:screening}
\end{center}
\end{table}

Table \ref{tab:screening} shows that the upgraded Algorithm~\ref{alg:IterativeCoefStreng} screens out a large percentage of the constraints in model \eqref{eq:BigM-OPF}, only retaining between 2\% and 26\% of them. The results in this table demonstrate that combining the tightening of the Big-Ms and the elimination of superfluous constraints leads to significant computational savings. For instance, the speedup factor for the three smallest systems ranges now between 1.5x and 15.1x. In addition,  \textbf{TS}(3) is able to solve six and four instances to optimality for systems IEEE-118 and IEEE-300, respectively, in less than 10 hours, and the average MIP-GAP is reduced to 0.01\% in these two largest power systems. Therefore, the computational performance of model \eqref{eq:BigM-OPF} has been drastically improved by combining the strengthening of the Big-Ms (making model \eqref{eq:BigM-OPF} tighter) and the removal of redundant constraints (making model \eqref{eq:BigM-OPF} more compact). Equally important, the elimination of superfluous constraints notably reduces the MIP solver's need for RAM memory, which decreases from around 100 GB in \textbf{T} to 32 GB in~\textbf{TS}.

We continue the numerical experiments by evaluating the impact of including the valid inequalities derived in Section \ref{sec:Valid} in model \eqref{eq:BigM-OPF}. The so obtained results are collated in Table \ref{tab:valid0}. In what follows, this approach is called \textbf{BN+V} for short. The results in Table \ref{tab:valid0} comprise, in order and following the previous notation, the average  number of constraints, the average values of LR-GAP and MIP-GAP, the number of instances solved to optimality, and the average speedup factor with respect to the \textbf{BN} approach of Table \ref{tab:bench}.  

\begin{table}[htbp]
\begin{center}
\caption{Tightening by valid inequalities (\textbf{BN+V}): Results}
\begin{tabular}{cccccc}
\hline
 &IEEE-RTS-24  &IEEE-57 &IEEE-RTS-73 &IEEE-118 &IEEE-300\\
\hline
\#CON &101.0\% &101.2\% &101.2\% &101.6\%	&101.5\%\\
LR-GAP &0.3374\%	&0.2038\%	&0.0001\%	&0.4784\%	&0.3192\%\\
MIP-GAP (\#OPT) &0.00\% (10) &0.00\% (10) &0.00\% (10) &0.03\% (1) &0.08\% (0)\\
Speedup factor &46.4x	&13.8x	&65.2x	&1.1x	&1.0x\\
\hline
\end{tabular}
\label{tab:valid0}
\end{center}
\end{table}

As can be seen, our valid inequalities only increase the total number of constraints by 1.0-1.6\%. However, the LR-GAP is significantly reduced compared to that obtained by \textbf{BN}. This effect is particularly noticeable for the IEEE-RTS-73 system with an average value of the linear relaxation gap equal to 0.0001\%, meaning that the linear relaxation of problem~\eqref{eq:BigM-OPF} with the proposed valid inequalities is very tight and its solution very close to the actual solution of that problem. Furthermore, the inclusion of the valid inequalities lead to average speedup factors that range between 13.8x and 65.2x for the three smallest systems. For the two largest systems, the time limit is reached in most instances, but the average MIP-GAP is reduced to 0.03\% and 0.08\%, respectively. 

Finally, we present similar simulation results for the setup in which the valid inequalities are also used to boost the tightening and screening power of Algorithm \ref{alg:IterativeCoefStreng} (that is, the valid inequalities are also included in problem \eqref{RelMIP}), leading to method \textbf{TS+V}. Table~\ref{tab:valid} provides the average number of constraints, the average values of LR-GAP and MIP-GAP, and the average speedup factor of \textbf{TS+V} with respect to \textbf{BN} in Table~\ref{tab:bench}. Since increasing the number of iterations of Algorithm~\ref{alg:IterativeCoefStreng} barely affects the performance of \textbf{TS+V}, all the data shown in Table \ref{tab:valid} correspond to $\kappa=1$. 

\begin{table}[htbp]
\begin{center}
\caption{Tightening and screening with valid inequalities (\textbf{TS+V}): Results}
\begin{tabular}{ccccccc}
\hline
 &IEEE-RTS-24  &IEEE-57 &IEEE-RTS-73 &IEEE-118 &IEEE-300\\
 \hline
\#CON &3.49\% &1.61\% &3.84\% &2.68\% &3.30\%\\
LR-GAP &0.3365\%	&0.1519\%	&0.0001\%	&0.2821\%	&0.1603\%\\
MIP-GAP (\#OPT) &0.00\% (10) &0.00\% (10) &0.00\% (10) &0.00\% (10) &0.00\% (10)\\
Speedup factor &706.8x	&35.3x &1470.0x	&23.1x	&8.5x\\
\hline
\end{tabular}
\label{tab:valid}
\end{center}
\end{table}

The comparison of the results in Tables~\ref{tab:screening}, \ref{tab:valid0} and \ref{tab:valid} yields the following observations. First, including the valid inequalities in Algorithm~\ref{alg:IterativeCoefStreng} strengthens the Big-Ms even further, which in turn remarkably reduces the linear relaxation gap and increases the number of constraints identified as redundant in model \eqref{eq:BigM-OPF}. Indeed, the total number of constraints eventually retained by \textbf{TS+V} ranges between 1.61\% and 3.84\% if compared with \textbf{BN}. Second, \textbf{TS+V} can solve the 10 instances to global optimality in less than 10 hours for the five test systems considered in these numerical experiments. In fact, the optimal solutions obtained by \textbf{TS+V} are the ones we use to compute the values of the linear relaxation gap throughout these simulations. Third, \textbf{TS+V} is able to achieve speedup factors between 8.5x and 1470.0x depending on the test system. All in all, \textbf{TS+V} features the best computational performance in terms of resolution time and MIP-GAP among all the methods tested so far.

Next, in Table~\ref{tab:comparison}  the results of \textbf{TS+V} are contrasted with those provided by state-of-the-art approximations available in the literature. In particular, we consider the CVaR-based and ALSO-X conservative approximations, both described in \cite{jiang2022}. Table~\ref{tab:comparison} provides the average cost increase in percentage with respect to the optimal cost and the speedup factor with regard to \textbf{BN} for the different methodologies compared. As expected, the CVaR-based approximation leads to conservative results and involves average cost increases that range between 0.49\% and 2.88\%. Interestingly enough, \textbf{TS+V} computes the optimal solution and involves a higher speedup factor than the CVaR-based approach for two of the five systems. Compared with the two ALSO-X approximations, \textbf{TS+V} obtains the global optimal solution in all cases with speedup factors that are still tantamount to those of these approximate methods.

\begin{table}[htbp]
\begin{center}
\caption{Comparison of the proposed \textbf{TS+V} approach and existing approximate methods}
\begin{tabular}{p{2.5cm}lccccc}
\hline
\multicolumn{2}{c}{} &IEEE-RTS-24  &IEEE-57 &IEEE-RTS-73 &IEEE-118 &IEEE-300\\
\hline
\multirow{4}{=}{Average cost increase}
&\textbf{TS+V} &0.00\% &0.00\%	&0.00\%	&0.00\%	&0.00\%\\
&CVaR &2.88\% &0.53\%	&1.71\%	&0.57\%	&0.49\%\\
&ALSO-X &0.80\%	&0.08\%	&0.41\%	&0.08\%	&0.05\%\\
&ALSO-X+ &0.53\%	&0.07\%	&0.12\%	&0.05\%	&0.04\%\\
\hline
\multirow{4}{=}{Speedup factor}
&\textbf{TS+V}  &706.8x	&35.3x &1470.0x	&23.1x	&8.5x\\
&CVaR &387.1x	&117.0x &265.1x	&4045.5x	&779.7x\\
&ALSO-X &18.0x	&1.2x &21.3x	&148.6x	&11.13x\\
&ALSO-X+ &7.8x	&0.7x &6.6x	&49.3x	&3.7x\\
\hline
\end{tabular}
\label{tab:comparison}
\end{center}
\end{table}

To conclude this study, we compare our solution approach with the one proposed in \cite{luedtke2014}. This author presents an exact approach to solve CCPs based on the addition of quantile cuts in a so-called \textit{lazy fashion} to a relaxation of the problem that only includes the deterministic constraints. In this way, at each node of the branching tree we seek for a violated constraint from a scenario $s$ with $y_s=0$. Using the coefficients of this violated constraint, a valid inequality is derived by means of the resolution of $|\mathcal{S}|$  linear separation subproblems. These strengthened cuts contain binary variables associated to several scenarios and are obtained by applying the \textit{star inequalities} of~\cite{atamturk2000}.

We have implemented the procedure including, at each iteration, both the initial or \textit{base} violated inequality and the strengthened one, since this approach delivered the best results. The comparison between our approach and Luedtke's (which we name BCD for short, from \emph{branch-and-cut decomposition}) is shown in Table \ref{tab:comparisonLuedtke}. As can be seen, a similar number of constraints are generated in the BCD approach (that includes the deterministic constraints and the cuts generated dynamically). Furthermore, our procedure clearly outperforms the BCD approach in terms of  computational time. One possible explanation, already pointed out in the paper, is that the generation of the cuts requires solving 1000 subproblems at each iteration. In the particular application presented in \cite{luedtke2014}, however, this potential bottleneck is overcome by leveraging the specific structure of the subproblem to develop a closed-form expression of its solution, thus avoiding the need to optimize over all the scenario sets. With that said, we remark that the times provided in Table~\ref{tab:comparisonLuedtke} under the acronym BCD-P correspond to a \emph{parallel} implementation of the BCD algorithm, in which only the time required to solve the most time-consuming subproblem out of the 1000 to be solved at each iteration is added to the final solution time reported. For completeness and because the \emph{parallel} implementation of BCD using off-the-shelf optimization solvers is not trivial at all, we also report in Table~\ref{tab:comparisonLuedtke} the speedup factor that is attained by a \emph{serial} implementation of BCD, termed  BCD-S in the table. As can be seen, this implementation is, however, far from being competitive in all cases.

\begin{table}[htbp]
\begin{center}
\caption{Comparison of the proposed \textbf{TS+V} approach and the BCD approach} \label{}
\begin{tabular}{llccccc}
\hline
\multicolumn{2}{c}{} &IEEE-RTS-24  &IEEE-57 &IEEE-RTS-73  &IEEE-118 &IEEE-300\\
\hline
\multirow{2}{*}{\#CONS}
&\textbf{TS+V} &3.49\% &1.61\% &3.84\% &2.68\% &3.30\%\\
&BCD &11.37\% &2.69\%	&3.58\%	&2.46\%	&1.53\%\\
\hline
\multirow{3}{*}{Speedup factor}
&\textbf{TS+V}  &706.8x	&35.3x &1470.0x &23.1x	&8.5x\\
&BCD-P &10.12x	&2.05x &69.00x	&9.24x	&4.74x\\
&BCD-S &0.12x	&0.04x &0.30x	&1.45x	&0.42x\\
\hline
\end{tabular}
\label{tab:comparisonLuedtke}
\end{center}
\end{table}

\section{Conclusions and future research}  \label{sec:conclusion}
In this paper, we propose a novel exact resolution technique for a MIP SAA-based reformulation of joint chance-constrained problems in the form of~\eqref{GenCC}. Our methodology includes a screening method to eliminate superfluous constraints based on an iterative procedure to repeatedly tighten the Big-Ms present in the MIP. These procedures are combined with the addition of inequalities that are valid provided that the technology matrix in the chance-constraint system exhibits the structure presented in~\eqref{MIPGenCC}. Said inequalities strengthen the  linear relaxation of the MIP SAA-based reformulation and allow for additional screening of constraints. The resultant model is thus compact and tight.

We have applied our methodology to solve the joint chance-constrained DC Optimal Power Flow. In the case study, we show that, in comparison with the benchmark model, our methodology provides remarkable results in terms of the linear relaxation bounds, the RAM memory needed to solve the instances, and the total computational resolution time. Specifically, our method \textbf{TS+V} solves to optimality all of the instances generated for the IEEE-RTS-118 and the IEEE-RTS-300 test systems, the majority of which are not solved within 10 hours of computational time using the benchmark approach. Furthermore, the average number of constraints eliminated from all instances with \textbf{TS+V} always exceeds a 95\% of them, and the lower bound is markedly increased by the inclusion of the valid inequalities, showing the outstanding results of the combination of the methods developed.

The comparison of our results with those provided by existing approximate and exact methods shows that our approach is computationally very competitive for small and medium-sized instances, always providing the best results in terms of cost. For the large instances addressed, while outperformed by the approximate methods in terms of computational time (as expected), our exact solution strategy not only provides a certificate of optimality, but also returns the optimal solution within the set time limit. A promising future research line consists in the development of a generalized set of valid inequalities that combine variables from pairs or subgroups of constraints $j\in \mathcal{J}$ in the chance-constrained system \eqref{GenCC_chance}.

\section*{Acknowledgements} 
This work was supported in part by the European Research Council (ERC) under the EU Horizon 2020 research and innovation program (grant agreement No. 755705), in part by the Spanish Ministry of Science and Innovation (AEI/10.13039/501100011033) through project PID2020-115460GB-I00, and in part by the Junta de Andalucía (JA) and the European Regional Development Fund (FEDER) through the research project P20\_00153. \'A. Porras is also financially supported by the Spanish Ministry of Science, Innovation and Universities through the University Teacher Training Program with fellowship number FPU19/03053. Finally, the authors thankfully acknowledge the computer resources, technical expertise, and assistance provided by the SCBI (Supercomputing and Bioinformatics) center of the University of M\'alaga.

\bibliographystyle{informs2014}
\bibliography{BibtexCC}

\begin{thebibliography}{75}
\providecommand{\natexlab}[1]{#1}
\providecommand{\url}[1]{\texttt{#1}}
\providecommand{\urlprefix}{URL }

\bibitem[{Abdi \protect\BIBand{} Fukasawa(2016)}]{abdi2016}
Abdi A, Fukasawa R (2016) On the mixing set with a knapsack constraint.
  \emph{Mathematical Programming} 157(1):191--217.

\bibitem[{Agarwal et~al.(1998)Agarwal, Aronov, Chan, \protect\BIBand{}
  Sharir}]{agarwal1998}
Agarwal PK, Aronov B, Chan TM, Sharir M (1998) On {{levels}} in
  {{arrangements}} of {{lines}}, {{segments}}, {{planes}}, and {{triangles}}.
  \emph{Discrete \& Computational Geometry} 19(3):315--331.

\bibitem[{Agarwal \protect\BIBand{} Matou{\v s}ek(1995)}]{agarwal1995}
Agarwal PK, Matou{\v s}ek J (1995) Dynamic half-space range reporting and its
  applications. \emph{Algorithmica} 13(4):325--345.

\bibitem[{Ahmed et~al.(2017)Ahmed, Luedtke, Song, \protect\BIBand{}
  Xie}]{ahmed2017}
Ahmed S, Luedtke J, Song Y, Xie W (2017) Nonanticipative duality, relaxations,
  and formulations for chance-constrained stochastic programs.
  \emph{Mathematical Programming} 162(1):51--81.

\bibitem[{Ahmed \protect\BIBand{} Xie(2018)}]{ahmed2018}
Ahmed S, Xie W (2018) Relaxations and approximations of chance constraints
  under finite distributions. \emph{Mathematical Programming} 170(1):43--65.

\bibitem[{Atamt{\"u}rk et~al.(2000)Atamt{\"u}rk, Nemhauser, \protect\BIBand{}
  Savelsbergh}]{atamturk2000}
Atamt{\"u}rk A, Nemhauser GL, Savelsbergh MWP (2000) Conflict graphs in solving
  integer programming problems. \emph{European Journal of Operational Research}
  121(1):40--55, ISSN 0377-2217,
  \urlprefix\url{http://dx.doi.org/10.1016/S0377-2217(99)00015-6}.

\bibitem[{{Ben-Tal} \protect\BIBand{} Nemirovski(2000)}]{ben-tal2000}
{Ben-Tal} A, Nemirovski A (2000) Robust solutions of linear programming
  problems contaminated with uncertain data. \emph{Mathematical programming}
  88(3):411--424.

\bibitem[{Bienstock et~al.(2014)Bienstock, Chertkov, \protect\BIBand{}
  Harnett}]{bienstock2014chance}
Bienstock D, Chertkov M, Harnett S (2014) Chance-constrained optimal power
  flow: Risk-aware network control under uncertainty. \emph{SIAM Review}
  56(3):461--495.

\bibitem[{Borobia(1994)}]{borobia1994}
Borobia A (1994) Mirror property for nonsingular mixed configurations of lines
  and points in {{$\mathbb{R}$3}}. \emph{Discrete \& Computational Geometry}
  11(3):311--320.

\bibitem[{Chan et~al.(1997)Chan, Snoeyink, \protect\BIBand{} Yap}]{chan1997}
Chan TM, Snoeyink J, Yap CK (1997) Primal dividing and dual pruning:
  {{Output-sensitive}} construction of four-dimensional polytopes and
  three-dimensional {{Voronoi}} diagrams. \emph{Discrete \& Computational
  Geometry} 18(4):433--454.

\bibitem[{Chazelle \protect\BIBand{} Edelsbrunner(1987)}]{chazelle1987}
Chazelle B, Edelsbrunner H (1987) An improved algorithm for constructing
  kth-order {{Voronoi}} diagrams. \emph{IEEE Transactions on Computers}
  100(11):1349--1354.

\bibitem[{Chazelle \protect\BIBand{} Preparata(1986)}]{chazelle1986}
Chazelle B, Preparata FP (1986) Halfspace range search: An algorithmic
  application of k-sets. \emph{Discrete \& Computational Geometry} 1(1):83--93.

\bibitem[{Cheema et~al.(2014)Cheema, Shen, Lin, \protect\BIBand{}
  Zhang}]{cheema2014}
Cheema MA, Shen Z, Lin X, Zhang W (2014) A unified framework for efficiently
  processing ranking related queries. \emph{EDBT}, 427--438.

\bibitem[{Chen et~al.(2021)Chen, Zhang, Hui, \protect\BIBand{}
  Song}]{chen2021time}
Chen G, Zhang H, Hui H, Song Y (2021) Time-efficient joint chance-constrained
  optimal power flow with a learning-based robust approximation. \emph{arXiv
  preprint arXiv:2112.09827} .

\bibitem[{Clarkson(1987)}]{clarkson1987}
Clarkson KL (1987) New applications of random sampling in computational
  geometry. \emph{Discrete \& Computational Geometry} 2(2):195--222.

\bibitem[{Conforti et~al.(2014)Conforti, Cornu{\'e}jols, \protect\BIBand{}
  Zambelli}]{conforti2014}
Conforti M, Cornu{\'e}jols G, Zambelli G (2014) \emph{Integer {{Programming}}},
  volume 271 ({Springer}).

\bibitem[{Dan{\'i}elsson et~al.(2008)Dan{\'i}elsson, Jorgensen, {de Vries},
  \protect\BIBand{} Yang}]{danielsson2008}
Dan{\'i}elsson J, Jorgensen BN, {de Vries} CG, Yang X (2008) Optimal portfolio
  allocation under the probabilistic {{VaR}} constraint and incentives for
  financial innovation. \emph{Annals of Finance} 4(3):345--367.

\bibitem[{Das et~al.(2007)Das, Gunopulos, Koudas, \protect\BIBand{}
  Sarkas}]{das2007}
Das G, Gunopulos D, Koudas N, Sarkas N (2007) Ad-hoc top-k query answering for
  data streams. \emph{Proceedings of the 33rd International Conference on
  {{Very}} Large Data Bases}, 183--194.

\bibitem[{Dentcheva et~al.(2000)Dentcheva, Pr{\'e}kopa, \protect\BIBand{}
  Ruszczynski}]{dentcheva2000}
Dentcheva D, Pr{\'e}kopa A, Ruszczynski A (2000) Concavity and efficient points
  of discrete distributions in probabilistic programming. \emph{Mathematical
  Programming} 89(1):55--77.

\bibitem[{Dey(1998)}]{dey1998}
Dey TK (1998) Improved bounds for planar k-sets and related problems.
  \emph{Discrete \& Computational Geometry} 19(3):373--382.

\bibitem[{Edelsbrunner \protect\BIBand{} Welzl(1985)}]{edelsbrunner1985}
Edelsbrunner H, Welzl E (1985) On the number of line separations of a finite
  set in the plane. \emph{Journal of Combinatorial Theory, Series A}
  38(1):15--29.

\bibitem[{Edelsbrunner \protect\BIBand{} Welzl(1986)}]{edelsbrunner1986}
Edelsbrunner H, Welzl E (1986) Constructing belts in two-dimensional
  arrangements with applications. \emph{SIAM Journal on Computing} 15(1).

\bibitem[{El{\c c}i et~al.(2018)El{\c c}i, Noyan, \protect\BIBand{}
  B{\"u}lb{\"u}l}]{elci2018a}
El{\c c}i {\"O}, Noyan N, B{\"u}lb{\"u}l K (2018) Chance-constrained stochastic
  programming under variable reliability levels with an application to
  humanitarian relief network design. \emph{Computers \& Operations Research}
  96:91--107.

\bibitem[{Esteban-P{\'e}rez \protect\BIBand{}
  Morales(2023)}]{esteban2023distributionally}
Esteban-P{\'e}rez A, Morales JM (2023) Distributionally robust optimal power
  flow with contextual information. \emph{European Journal of Operational
  Research} 306(3):1047--1058.

\bibitem[{Frank et~al.(2012)Frank, Steponavice, \protect\BIBand{}
  Rebennack}]{frank2012optimal}
Frank S, Steponavice I, Rebennack S (2012) Optimal power flow: {A}
  bibliographic survey {I}. \emph{Energy Systems} 3(3):221--258.

\bibitem[{G{\"u}nl{\"u}k \protect\BIBand{} Pochet(2001)}]{gunluk2001}
G{\"u}nl{\"u}k O, Pochet Y (2001) Mixing mixed-integer inequalities.
  \emph{Mathematical Programming} 90(3):429--457.

\bibitem[{{Gurobi Optimization, LLC}(2022)}]{gurobi}
{Gurobi Optimization, LLC} (2022) {Gurobi Optimizer Reference Manual}.
  \urlprefix\url{https://www.gurobi.com}.

\bibitem[{Halperin \protect\BIBand{} Sharir(2017)}]{halperin2017}
Halperin D, Sharir M (2017) Arrangements. \emph{Handbook of {{Discrete}} and
  {{Computational Geometry}}}, 723--762 ({Chapman and Hall/CRC}), ISBN
  978-1-315-11960-1,
  \urlprefix\url{http://dx.doi.org/10.1201/9781315119601-28}.

\bibitem[{Henrion(2007)}]{henrion2007}
Henrion R (2007) Structural properties of linear probabilistic constraints.
  \emph{Optimization} 56(4):425--440.

\bibitem[{Henrion \protect\BIBand{} Strugarek(2011)}]{henrion2011}
Henrion R, Strugarek C (2011) Convexity of chance constraints with dependent
  random variables: The use of copulae. \emph{{Stochastic Optimization Methods
  in Finance and Energy}}, 427--439 ({Springer}).

\bibitem[{Hong et~al.(2011)Hong, Yang, \protect\BIBand{} Zhang}]{hong2011}
Hong LJ, Yang Y, Zhang L (2011) Sequential convex approximations to joint
  chance constrained programs: {{A Monte Carlo}} approach. \emph{Operations
  Research} 59(3):617--630.

\bibitem[{Hou \protect\BIBand{} Roald(2020)}]{LineTunning}
Hou AM, Roald LA (2020) Chance constraint tuning for optimal power flow.
  \emph{2020 International Conference on Probabilistic Methods Applied to Power
  Systems (PMAPS)}, 1--6,
  \urlprefix\url{http://dx.doi.org/10.1109/PMAPS47429.2020.9183552}.

\bibitem[{Jiang \protect\BIBand{} Xie(2022)}]{jiang2022}
Jiang N, Xie W (2022) {ALSO-X} and {ALSO-X+}: Better convex approximations for
  chance constrained programs. \emph{Operations Research} 70(6):3581--3600.

\bibitem[{Küçükyavuz \protect\BIBand{} Jiang(2022)}]{kucukyavuz2022}
Küçükyavuz S, Jiang R (2022) Chance-constrained optimization under limited
  distributional information: {A} review of reformulations based on sampling
  and distributional robustness. \emph{EURO Journal on Computational
  Optimization} 10:100030, ISSN 2192-4406,
  \urlprefix\url{http://dx.doi.org/https://doi.org/10.1016/j.ejco.2022.100030}.

\bibitem[{Lagoa et~al.(2005)Lagoa, Li, \protect\BIBand{} Sznaier}]{lagoa2005}
Lagoa CM, Li X, Sznaier M (2005) Probabilistically constrained linear programs
  and risk-adjusted controller design. \emph{SIAM Journal on Optimization}
  15(3):938--951.

\bibitem[{Lejeune \protect\BIBand{} Dehghanian(2020)}]{lejeune2020optimal}
Lejeune MA, Dehghanian P (2020) Optimal power flow models with probabilistic
  guarantees: A boolean approach. \emph{IEEE Transactions on Power Systems}
  35(6):4932--4935.

\bibitem[{Luedtke(2014)}]{luedtke2014}
Luedtke J (2014) A branch-and-cut decomposition algorithm for solving
  chance-constrained mathematical programs with finite support.
  \emph{Mathematical Programming} 146(1-2):219--244, ISSN 0025-5610, 1436-4646,
  \urlprefix\url{http://dx.doi.org/10.1007/s10107-013-0684-6}.

\bibitem[{Luedtke \protect\BIBand{} Ahmed(2008)}]{luedtke2008}
Luedtke J, Ahmed S (2008) A sample approximation approach for optimization with
  probabilistic constraints. \emph{SIAM Journal on Optimization}
  19(2):674--699.

\bibitem[{Luedtke et~al.(2010)Luedtke, Ahmed, \protect\BIBand{}
  Nemhauser}]{luedtke2010}
Luedtke J, Ahmed S, Nemhauser GL (2010) An integer programming approach for
  linear programs with probabilistic constraints. \emph{Mathematical
  programming} 122(2):247--272.

\bibitem[{Miller \protect\BIBand{} Wagner(1965)}]{miller1965chance}
Miller BL, Wagner HM (1965) Chance constrained programming with joint
  constraints. \emph{Operations Research} 13(6):930--945.

\bibitem[{Nair \protect\BIBand{} {Miller-Hooks}(2011)}]{nair2011}
Nair R, {Miller-Hooks} E (2011) Fleet management for vehicle sharing
  operations. \emph{Transportation Science} 45(4):524--540.

\bibitem[{Najjarbashi \protect\BIBand{} Lim(2020)}]{najjarbashi2020}
Najjarbashi A, Lim GJ (2020) A decomposition algorithm for the two-stage
  chance-constrained operating room scheduling problem. \emph{IEEE Access}
  8:80160--80172.

\bibitem[{Natarajan et~al.(2008)Natarajan, Pachamanova, \protect\BIBand{}
  Sim}]{natarajan2008}
Natarajan K, Pachamanova D, Sim M (2008) Incorporating asymmetric
  distributional information in robust value-at-risk optimization.
  \emph{Management Science} 54(3):573--585.

\bibitem[{Nemirovski \protect\BIBand{} Shapiro(2006)}]{nemirovski2006}
Nemirovski A, Shapiro A (2006) Scenario approximations of chance constraints.
  \emph{Probabilistic and randomized methods for design under uncertainty}
  3--47.

\bibitem[{Nemirovski \protect\BIBand{} Shapiro(2007)}]{nemirovski2007}
Nemirovski A, Shapiro A (2007) Convex approximations of chance constrained
  programs. \emph{SIAM Journal on Optimization} 17(4):969--996.

\bibitem[{Nivasch(2008)}]{nivasch2008}
Nivasch G (2008) An improved, simple construction of many halving edges.
  \emph{Contemporary Mathematics} 453:299--306.

\bibitem[{OASYS(2023)}]{TCSAAJCCOPF2023}
OASYS (2023) {Data and Code for a Tight and Compact Model of the {SAA}-based
  Joint Chance-constrained {OPF}}. \emph{GitHub repository
  (https://github.com/groupoasys/TC{\_}SAA{\_}JCC-OPF)}
  \urlprefix\url{https://github.com/groupoasys/TC_SAA_JCC-OPF}.

\bibitem[{Peña-Ordieres et~al.(2020)Peña-Ordieres, Luedtke, \protect\BIBand{}
  Wächter}]{pena-ordieres2020}
Peña-Ordieres A, Luedtke JR, Wächter A (2020) Solving chance-constrained
  problems via a smooth sample-based nonlinear approximation. \emph{SIAM
  Journal on Optimization} 30(3):2221--2250,
  \urlprefix\url{http://dx.doi.org/10.1137/19M1261985}.

\bibitem[{Peña-Ordieres et~al.(2021)Peña-Ordieres, Molzahn, Roald,
  \protect\BIBand{} Wächter}]{LineAlejandra}
Peña-Ordieres A, Molzahn DK, Roald LA, Wächter A (2021) {DC} optimal power
  flow with joint chance constraints. \emph{IEEE Transactions on Power Systems}
  36(1):147--158, \urlprefix\url{http://dx.doi.org/10.1109/TPWRS.2020.3004023}.

\bibitem[{{Power Grid Lib}(2022)}]{pglib}
{Power Grid Lib} (2022) GitHub repository, available at:
  (https://github.com/power-grid-lib/pglib-opf).

\bibitem[{Pr{\'e}kopa(2003)}]{prekopa2003}
Pr{\'e}kopa A (2003) Probabilistic programming. \emph{Handbooks in Operations
  Research and Management Science} 10:267--351.

\bibitem[{Preparata \protect\BIBand{} Shamos(2012)}]{preparata2012}
Preparata FP, Shamos MI (2012) \emph{Computational Geometry: An Introduction}
  ({Springer Science \& Business Media}).

\bibitem[{Qiu et~al.(2014)Qiu, Ahmed, Dey, \protect\BIBand{} Wolsey}]{qiu2014}
Qiu F, Ahmed S, Dey SS, Wolsey LA (2014) Covering linear programming with
  violations. \emph{INFORMS Journal on Computing} 26(3):531--546,
  \urlprefix\url{http://dx.doi.org/10.1287/ijoc.2013.0582}.

\bibitem[{Ramos(1999)}]{ramos1999}
Ramos EA (1999) On range reporting, ray shooting and k-level construction.
  \emph{Proceedings of the Fifteenth Annual Symposium on {{Computational
  Geometry}}}, 390--399.

\bibitem[{Roald \protect\BIBand{} Andersson(2018)}]{LineGoran}
Roald L, Andersson G (2018) Chance-constrained {AC} optimal power flow:
  {R}eformulations and efficient algorithms. \emph{IEEE Transactions on Power
  Systems} 33(3):2906--2918,
  \urlprefix\url{http://dx.doi.org/10.1109/TPWRS.2017.2745410}.

\bibitem[{Rockafellar \protect\BIBand{} Uryasev(2000)}]{rockafellar2000}
Rockafellar RT, Uryasev S (2000) Optimization of conditional value-at-risk.
  \emph{The Journal of Risk} 2(3):21--41, ISSN 14651211,
  \urlprefix\url{http://dx.doi.org/10.21314/JOR.2000.038}.

\bibitem[{Roos \protect\BIBand{} Widmayer(1994)}]{roos1994}
Roos T, Widmayer P (1994) K-{{Violation}} linear programming. \emph{Information
  Processing Letters} 52(2):109--114, ISSN 0020-0190,
  \urlprefix\url{http://dx.doi.org/10.1016/0020-0190(94)00134-0}.

\bibitem[{Shapiro(2003)}]{shapiro2003}
Shapiro A (2003) Monte {{Carlo}} sampling methods. \emph{Handbooks in
  operations research and management science} 10:353--425.

\bibitem[{Sharir(2011)}]{sharir2011}
Sharir M (2011) An improved bound for k-sets in four dimensions.
  \emph{Combinatorics, Probability and Computing} 20(1):119--129.

\bibitem[{Sharir et~al.(2000)Sharir, Smorodinsky, \protect\BIBand{}
  Tardos}]{sharir2000}
Sharir M, Smorodinsky S, Tardos G (2000) An improved bound for k-sets in three
  dimensions. \emph{Proceedings of the Sixteenth Annual Symposium on
  {{Computational}} Geometry}, 43--49.

\bibitem[{Shaw et~al.(2016)Shaw, Irfan, Shankar, \protect\BIBand{}
  Yadav}]{shaw2016}
Shaw K, Irfan M, Shankar R, Yadav SS (2016) Low carbon chance constrained
  supply chain network design problem: A {{Benders}} decomposition based
  approach. \emph{Computers \& Industrial Engineering} 98:483--497.

\bibitem[{Shen et~al.(2012)Shen, Cheema, \protect\BIBand{} Lin}]{shen2012}
Shen Z, Cheema MA, Lin X (2012) Loyalty-based selection: {{Retrieving}} objects
  that persistently satisfy criteria. \emph{Proceedings of the 21st {{ACM}}
  International Conference on {{Information}} and Knowledge Management},
  2189--2193.

\bibitem[{Song et~al.(2014)Song, Luedtke, \protect\BIBand{} K{\"u}{\c
  c}{\"u}kyavuz}]{song2014}
Song Y, Luedtke JR, K{\"u}{\c c}{\"u}kyavuz S (2014) Chance-constrained binary
  packing problems. \emph{INFORMS Journal on Computing} 26(4):735--747.

\bibitem[{Sun et~al.(2014)Sun, Xu, \protect\BIBand{} Wang}]{sun2014}
Sun H, Xu H, Wang Y (2014) Asymptotic analysis of sample average approximation
  for stochastic optimization problems with joint chance constraints via
  conditional value at risk and difference of convex functions. \emph{Journal
  of Optimization Theory and Applications} 161(1):257--284, ISSN 0022-3239,
  1573-2878, \urlprefix\url{http://dx.doi.org/10.1007/s10957-012-0127-1}.

\bibitem[{Taleizadeh et~al.(2012)Taleizadeh, Niaki, \protect\BIBand{}
  Makui}]{taleizadeh2012}
Taleizadeh AA, Niaki STA, Makui A (2012) Multiproduct multiple-buyer
  single-vendor supply chain problem with stochastic demand, variable
  lead-time, and multi-chance constraint. \emph{Expert Systems with
  Applications} 39(5):5338--5348.

\bibitem[{Tanner et~al.(2008)Tanner, Sattenspiel, \protect\BIBand{}
  Ntaimo}]{tanner2008}
Tanner MW, Sattenspiel L, Ntaimo L (2008) Finding optimal vaccination
  strategies under parameter uncertainty using stochastic programming.
  \emph{Mathematical Biosciences} 215(2):144--151, ISSN 0025-5564,
  \urlprefix\url{http://dx.doi.org/10.1016/j.mbs.2008.07.006}.

\bibitem[{Tayur et~al.(1995)Tayur, Thomas, \protect\BIBand{}
  Natraj}]{tayur1995}
Tayur SR, Thomas RR, Natraj NR (1995) An algebraic geometry algorithm for
  scheduling in presence of setups and correlated demands. \emph{Mathematical
  Programming} 69(1):369--401.

\bibitem[{T{\'o}th et~al.(2017)T{\'o}th, O'Rourke, \protect\BIBand{}
  Goodman}]{toth2017}
T{\'o}th CD, O'Rourke J, Goodman JE (2017) \emph{{Handbook of Discrete and
  Computational Geometry}} ({CRC press}).

\bibitem[{T{\'o}th(2000)}]{toth2000}
T{\'o}th G (2000) Point sets with many k-sets. \emph{Proceedings of the
  Sixteenth Annual Symposium on {{Computational}} Geometry}, 37--42.

\bibitem[{Van~Ackooij et~al.(2011)Van~Ackooij, Zorgati, Henrion,
  \protect\BIBand{} M{\"o}ller}]{vanackooij2011}
Van~Ackooij W, Zorgati R, Henrion R, M{\"o}ller A (2011) Chance constrained
  programming and its applications to energy management. \emph{Stochastic
  Optimization-Seeing the Optimal for the Uncertain} 291--320.

\bibitem[{Vrakopoulou et~al.(2013)Vrakopoulou, Margellos, Lygeros,
  \protect\BIBand{} Andersson}]{vrakopoulou2013}
Vrakopoulou M, Margellos K, Lygeros J, Andersson G (2013) {Probabilistic
  Guarantees for the {{N-1}} Security of Systems with Wind Power Generation}.
  \emph{{Reliability and Risk Evaluation of Wind Integrated Power Systems}},
  59--73 ({Springer}).

\bibitem[{Xie \protect\BIBand{} Ahmed(2016)}]{xie2016}
Xie W, Ahmed S (2016) On the quantile cut closure of chance-constrained
  problems. \emph{International {{Conference}} on {{Integer Programming}} and
  {{Combinatorial Optimization}}}, 398--409 ({Springer}).

\bibitem[{Xie \protect\BIBand{} Ahmed(2018)}]{xie2018}
Xie W, Ahmed S (2018) On quantile cuts and their closure for chance constrained
  optimization problems. \emph{Mathematical Programming} 172(1):621--646.

\bibitem[{Yu et~al.(2012)Yu, Agarwal, \protect\BIBand{} Yang}]{yu2012}
Yu A, Agarwal PK, Yang J (2012) Processing a large number of continuous
  preference top-k queries. \emph{Proceedings of the 2012 {{ACM SIGMOD}}
  International Conference on Management of Data}, 397--408.

\bibitem[{Zhang et~al.(2015)Zhang, Shen, \protect\BIBand{}
  Mathieu}]{zhang2015data}
Zhang Y, Shen S, Mathieu JL (2015) Data-driven optimization approaches for
  optimal power flow with uncertain reserves from load control. \emph{2015
  American Control Conference (ACC)}, 3013--3018 (IEEE).

\end{thebibliography}

\newpage
\setcounter{page}{1}





\appendix

\renewcommand{\thesection}{Appendix \Alph{section}}
\renewcommand{\thetheorem}{A.\arabic{theorem}}
\section{Relationship with \emph{quantile cuts}} \label{sec:anexoQuantileCuts}
\thispagestyle{empty}

Consider problem  \eqref{GenCC}, with feasible region $\mathcal{FR}$, assume $X$ is a non-empty and compact set given by  deterministic constraints 
and $a_j(\omega) = a^{0}_j + \Omega_j(\omega) \hat{a}_j$, $j \in \mathcal{J}$. For a given scenario $s\in \mathcal{S}$, define $\mathcal{X}_s$ as the set of points that satisfy the constraints for that scenario, i.e., 
\begin{equation*}
    \mathcal{X}_s := \left\{ x\in \mathbb{R}^{|\mathcal{I}|} : \Omega_{js} \hat{a}_j^{\top} x \le b_{js} - x^{\top} a_j^{0},\quad \forall j\in\mathcal{J}  \right\}.
\end{equation*}
To characterize the feasible region $\mathcal{FR}$, let $\mathcal{S}^F:= \left\{ S'\subseteq \mathcal{S}: |S'| \ge |\mathcal{S}|-p \right\}$ be the collection of all feasible subsets of scenarios. Then $\mathcal{FR}$ can be represented as:
\begin{equation*}
    \mathcal{FR} := \bigcup_{S'\in \mathcal{S}^F} \left[ X \bigcap_{s\in S'} \mathcal{X}_s \right].
\end{equation*}

In \cite{xie2018} it is proved that $ \mathcal{FR}$ can also be written in a conjunctive normal form. To do so, we introduce the definition of a \emph{minimal partial covering subset}.


\begin{definition}
A set $S''\subseteq \mathcal{S}$ is a \emph{partial covering subset} if it intersects with all the feasible scenario subsets in $\mathcal{S}^F$, i.e., for any $S'\in \mathcal{S}^F$ we have $S' \cap S'' \neq \emptyset$. Also, a set $S''$ is a \emph{minimal} partial covering subset if there does not exist another partial covering subset $\hat{S} \subseteq \mathcal{S}$ such that $\hat{S}\subsetneq S'$. Let $\mathcal{S}^C$ denote the collection of all minimal partial covering subsets. When the samples are independent and identically distributed, like in our case, $\mathcal{S}^C:= \left\{ S''\subseteq \mathcal{S}: |S''| = p+1 \right\}$.
\end{definition}

\begin{proposition}[{\cite[Proposition 3]{xie2018}}]
\begin{equation*}
 \mathcal{FR} =  \bigcap_{S''\in \mathcal{S}^C} \left[ \bigcup_{s\in S''} \left( X\bigcap \mathcal{X}_s \right) \right].
\end{equation*}
\end{proposition}

Consider now the relaxation of problem  \eqref{GenCC} associated to a given row $j\in \mathcal{J}$:
\begin{subequations} \label{GenCCSingle}
\begin{align}
 \min_{x} \quad &  f(x) \label{GenCCSingle_FO}\\
\text{s.t.}  \quad & \mathbb{P}\left\{ a_{j}(\omega)^{\top}x \le b_j(\omega) \right\} \ge 1-\epsilon \label{GenCCSingle_chance}\\
& x \in X \label{GenCCSingle_xinX} 
\end{align}
\end{subequations}
\noindent and its feasible region $\mathcal{FR}_j$. Clearly, it holds that $\mathcal{FR}\subseteq\bigcap_{j\in \mathcal{J}} \mathcal{FR}_j$. This is equivalent to considering the problem without joint constraints. In a similar manner, for a given scenario $s\in \mathcal{S}$ and row $j\in \mathcal{J}$, define $\mathcal{X}_{js}$ as the set of points that satisfy the scenario for that row:
\begin{equation*}
    \mathcal{X}_{js} := \left\{ x\in \mathbb{R}^{|\mathcal{I}|} : \Omega_{js} \hat{a}_j^{\top} x \le b_{js} - x^{\top} a_j^{0}  \right\},
\end{equation*}

\noindent with $\mathcal{X}_s=\bigcap_{j\in \mathcal{J}} \mathcal{X}_{js}$. Making use of the previous results, we have 
\begin{equation*}
    \mathcal{FR}_j := \bigcup_{S'\in \mathcal{S}^F} \left[ X \bigcap_{s\in S'} \mathcal{X}_{js} \right] = \bigcap_{S''\in \mathcal{S}^C} \left[ \bigcup_{s\in S''} \left( X\bigcap \mathcal{X}_{js} \right) \right].
\end{equation*}

Using this characterization of a relaxation of the feasible region, we are able to prove the following:
\begin{proposition}\label{prop:UequalFR}
Let $\mathcal{U}_j:= \left\{ x\in X: U_j^{p+1}(\hat{a}_j^{\top} x) \le - x^{\top} a_j^{0}\right\}$. Then for $j\in \mathcal{J}$, $\mathcal{U}_j = \mathcal{FR}_j$.
\end{proposition}

\proof
$[\mathcal{U}_j\subseteq \mathcal{FR}_j]$ Assume $x\in \mathcal{U}_j$. Then $x\in X$ satisfies the chance constraint for (at least) $|\mathcal{S}|-p$ scenarios. But then, for any minimal partial covering $S''\in\mathcal{S}^C$, since $|S''|=p+1$, we must have $x \in \bigcup_{s\in S''} (X \bigcap \mathcal{X}_{js})$. Given that this occurs for any $S''$, then $x$ belongs to the intersection of $S''$, and hence $x\in \mathcal{FR}_j$.

$[\mathcal{FR}_j\subseteq \mathcal{U}_j]$ This proof is straightforward considering the validity of the proof of Proposition \ref{prop:VVII_versionpiecewise}, but assuming $x\in \mathcal{FR}_j$ instead of $x\in \mathcal{FR}$. Thus, if $x\in \mathcal{FR}_j$, then $x$ satisfies the single chance constraint $j$ by definition. Since $\mathcal{U}_j^{p+1}(\hat{a}_j^{\top} x) \le - x^{\top} a_j^{0}$ is a valid  inequality for problem \eqref{GenCCSingle} associated to $j$, $x$ must satisfy it as well.
\endproof

Now, we are going to establish the relationship between our valid inequalities and the quantile closure described in \cite{xie2018}. Let us first introduce the necessary definitions to relate the quantile closure with $\mathrm{conv}(\mathcal{FR})$.
\begin{definition}
 Given $\alpha\in \mathbb{R}^{|\mathcal{I}|}$, let $\{\gamma^{\alpha}_s(X)\}$ be the optimal values of 
 \begin{equation*}
     \gamma^{\alpha}_s(X) := \min \left\{ \alpha^{\top}x: x\in X \cap \mathcal{X}_s\right\} \quad \forall s\in \mathcal{S}.
 \end{equation*}
 The $(1-\epsilon)$--\emph{quantile} (or simply \emph{quantile}) of $\{\gamma^{\alpha}_s(X)\}_{s\in S}$ is denoted by $\gamma^{\alpha}_q(X)$ and is given by
 \begin{equation*}
     \gamma^{\alpha}_q(X) := \min_{S'\in \mathcal{S}^F} \max_{s\in S'} \gamma^{\alpha}_s(X)
 \end{equation*}
 \noindent and the associated \emph{quantile cut} is 
 \begin{equation*}
     \alpha^{\top}x \ge \gamma^{\alpha}_q(X).
 \end{equation*}
 \end{definition}
 Note that the definition of quantile cuts depends on $X$, a set of deterministic constraints that we may change if we apply cuts successively. 

\begin{definition}
 The \emph{first quantile closure} of $X$ is defined as
 \begin{equation*}
     X^1 := \bigcap_{\alpha\in\mathbb{R}^{|\mathcal{I}|}} \left\{ x\in \mathbb{R}^{|\mathcal{I}|}: \alpha^{\top}x \ge \gamma_q^{\alpha}(X) \right\}.
 \end{equation*}
 Inductively, we define the \emph{$r$th round quantile closure} $X^r$ as:
 \begin{equation*}
     X^r := \bigcap_{\alpha\in\mathbb{R}^{|\mathcal{I}|}} \left\{ x\in \mathbb{R}^{|\mathcal{I}|}: \alpha^{\top}x \ge \gamma_q^{\alpha}(X^{r-1}) \right\} \quad r\ge 2.
 \end{equation*}
\end{definition}
\begin{theorem}[{\cite[Theorem 1]{xie2018}}]
\begin{equation} \label{eq:closureconjunctiveform}
    X^1 = \bigcap_{S''\in \mathcal{S}^C} \mathrm{conv} \left[ \bigcup_{s\in S''} \left( X \bigcap \mathcal{X}_s\right) \right],
\end{equation}
\noindent and for each $r\ge 2$, $r\in \mathbb{Z}_{++}$,
\begin{equation}  \label{eq:rclosureconjunctiveform}
     X^r = \bigcap_{S''\in \mathcal{S}^C} \mathrm{conv} \left[ \bigcup_{s\in S''} \left( X^{r-1} \bigcap \mathcal{X}_s\right) \right],
\end{equation}
\noindent where $X^0 = X$.
\end{theorem}
\begin{theorem}[{\cite[Theorem 2]{xie2018}}]
 The set sequence $\{X^r\}$ converges to $\mathrm{conv}(\mathcal{FR})$ with respect to the Hausdorff distance; i.e., $\widebar{X}:=\lim_{r\rightarrow\infty} X^r=\mathrm{conv}(\mathcal{FR})$.
\end{theorem}

If, like previously, we consider these definitions applied to the relaxation of problem \eqref{GenCC} that results from considering a single row $j\in \mathcal{J}$, then we can establish the following results:
\begin{definition}
 Given $\alpha\in \mathbb{R}^{|\mathcal{I}|}$, let $\{\gamma^{\alpha}_{js}(X)\}$ be the optimal values of 
 \begin{equation*}
     \gamma^{\alpha}_{js}(X) = \min \left\{ \alpha^{\top}x: x\in X \cap \mathcal{X}_{js}\right\} \quad \forall s\in \mathcal{S}.
 \end{equation*}
The $(1-\epsilon)$--\emph{quantile} (or simply \emph{quantile}) of $\{\gamma^{\alpha}_{js}(X)\}_{s\in S}$ is denoted by $\gamma^{\alpha}_{jq}(X)$ and is given by
 \begin{equation*}
     \gamma^{\alpha}_{jq}(X) := \min_{S'\in \mathcal{S}^F} \max_{s\in S'} \gamma^{\alpha}_{js}(X),
 \end{equation*}
 \noindent and the associated \emph{quantile cut} is 
 \begin{equation*}
     \alpha^{\top}x \ge \gamma^{\alpha}_{jq}(X).
 \end{equation*}
\end{definition}
\begin{definition}
 The \emph{first quantile closure} of $X$ associated to $j\in\mathcal{J}$ is defined as
 \begin{equation*}
     X^1_j := \bigcap_{\alpha\in\mathbb{R}^{|\mathcal{I}|}} \left\{ x\in \mathbb{R}^{|\mathcal{I}|}: \alpha^{\top}x \ge \gamma_{jq}^{\alpha}(X) \right\}.
 \end{equation*}
 Inductively, we define the \emph{$r$th round quantile closure} $X^r_j$ associated to $j\in\mathcal{J}$ as:
 \begin{equation*}
     X^r_j := \bigcap_{\alpha\in\mathbb{R}^{|\mathcal{I}|}} \left\{ x\in \mathbb{R}^{|\mathcal{I}|}: \alpha^{\top}x \ge \gamma_{jq}^{\alpha}(X^{r-1}) \right\} \quad r\ge 2.
 \end{equation*}
\end{definition}
\begin{theorem}[{\cite[Theorem 1]{xie2018}}]
\begin{equation} \label{eq:closureconjunctiveformj}
    X^1_j = \bigcap_{S''\in \mathcal{S}^C} \mathrm{conv} \left[ \bigcup_{s\in S''} \left( X \bigcap \mathcal{X}_{js}\right) \right],
\end{equation}
\noindent and for each $r\ge 2$, $r\in \mathbb{Z}_{++}$,
\begin{equation}  \label{eq:rclosureconjunctiveformj}
     X^r_j = \bigcap_{S''\in \mathcal{S}^C} \mathrm{conv} \left[ \bigcup_{s\in S''} \left( X^{r-1}_j \bigcap \mathcal{X}_{js}\right) \right],
\end{equation}
\noindent where $X^0_j = X$.
\end{theorem}
\begin{theorem}[{\cite[Theorem 2]{xie2018}}]
 The set sequence $\{X^r_j\}$ converges to $\mathrm{conv}(\mathcal{FR}_j)$ with respect to the Hausdorff distance; i.e., $\widebar{X}_j:=\lim_{r\rightarrow\infty} X^r_j=\mathrm{conv}(\mathcal{FR}_j)$.
\end{theorem}

\begin{corollary} \label{cor:clausuramedianteVVII}
\begin{equation}
    \mathcal{U}_j \subseteq \mathrm{conv} (\mathcal{U}_j) = \widebar{X}_j
\end{equation}
\end{corollary}

This result is particularly relevant because, as it is proven in the same work, the separation over the first quantile closure is, in general, NP-hard. Therefore, Corollary \ref{cor:clausuramedianteVVII} gives a way of deriving $\mathrm{conv}(\mathcal{FR}_j)$ through the explicit set of valid inequalities  \eqref{eq:VVII_versionupperhull}.

\section{Complexity and extension to the multidimensional case} \label{sec:anexoExtensionMultidim}

\renewcommand{\thetheorem}{B.\arabic{theorem}}
Suppose there exist matrices $A_j \in \mathbb{R}^{|\mathcal{I}| \times q_j }, q_j \leq |\mathcal{I}|$ (this includes the identity), $j \in \mathcal{J}$, such that $a_j(\omega)$ in the chance constraint can be written as $a_j(\omega) = a^{0}_j +  A_j \Omega_j(\omega)$, $j \in \mathcal{J}$, with $a_j^{0} \in \mathbb{R}^{|\mathcal{I}|}$ and $\Omega_j(\omega)$ being a $q_j$-dimensional vector of real-valued functions whose domains include the support of $\omega$.  Note that the case described in Subsection~\ref{sec:Valid} corresponds to $q_j = 1, \forall j \in \mathcal{J}$.

Now, define $z_j:=A_j^{\top} x$, with $z_j^d$ and $z_j^u$ being vectors containing the component-wise lower and upper bounds on $z$ induced by the feasibility set \eqref{MIPGenCC_xinX}. Then, we can define the sets of hyperplanes $\Pi_j$ for $j\in \mathcal{J}$
\begin{equation} \label{eq:setplanesLj}
    \Pi_j := \left\{ f_{js}(z_j) = \Omega_{js}^{\top} z_j - b_{js}, \quad \text{for } z_{j} \in B(z_j^d, z_j^u), \forall s \in \mathcal{S} \right\}, 
\end{equation} 
where $B(z_j^d, z_j^u):= \bigtimes\limits_{i=1}^{q_j}[z_{ji}^d, z_{ji}^u]$.

The $k$-lower envelope in higher dimension is defined as follows:
\begin{definition}
Given a set of $q$-dimensional hyperplanes $\Pi$, the lower (resp.\ upper) score of a point $t$ is the number of hyperplanes that lie strictly below (resp.\ strictly above) $t$. The $k$-lower (resp.\ $k$-upper) envelope is the closure of the set of points in $\Pi$ that have lower (resp.\ upper) score equal to $k-1$. The $k$-lower envelope is also known as the $k$-level.
\end{definition}

There is an extensive literature in the field of geometry dealing with the development of $k$-envelopes, since they find applications in seemingly unrelated problems like designing data structures (\cite{chazelle1986, clarkson1987, agarwal1995}), processing ranking queries (\cite{yu2012, das2007, shen2012, cheema2014}) or solving hyperplane partitioning problems (\cite{borobia1994, toth2000}). By geometric transformations, the time bounds found for $k$-envelopes also apply to closely related structures (and vice versa) like \emph{$k$-order Voronoi diagrams} (\cite{chan1997,chazelle1987,ramos1999}), or $k$-belts and $k$-sets: 

\begin{definition}
The \emph{$k$-belt} of a set of hyperplanes $\Pi$, for $0\le k\le \ceil*{\frac{|\Pi|}{2}}$, is the closure of the set of points in the space that have lower and upper score greater than or equal to $k$. 
\end{definition}
\begin{definition}
 Given a fixed set $T$ of points in $\mathbb{R}^d$, a \emph{$k$-set} ($k \le |T|$) of $T$ is any subset of $T$ of size $k$ which can be separated from the complement of $T$ by a hyperplane. The $k$-set problem seeks for the maximum possible number of $k$-sets. 
\end{definition}

\begin{figure}
\centering
\includegraphics[width=0.7\textwidth]{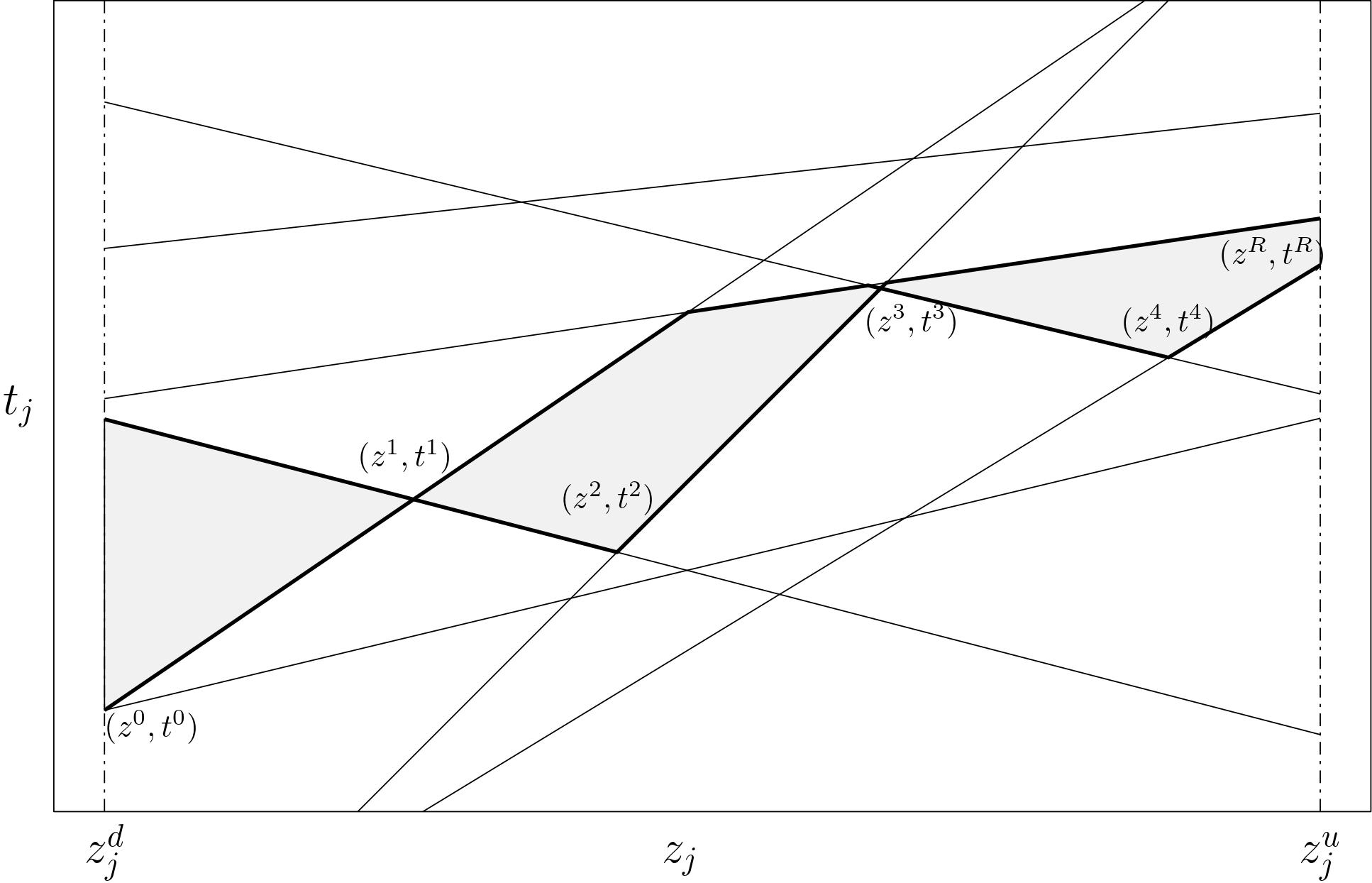} 
\caption{In grey, the 4-belt of a set $\mathcal{L}_j$ of 8 lines in the plane}\label{fig:kbelt} 
\end{figure}

Figure \ref{fig:kbelt} shows (in light grey) the $4$-belt of a set of 8 lines. Obviously, the boundary of the $k$-belt is the union of the $k$-upper and $k$-lower envelopes. Less obvious is the fact that $k$-envelopes and $k$-sets are closely related. In fact, using duality one can map a set of hyperplanes in $\mathbb{R}^d$ (the primal space) into a set of points in $\mathbb{R}^d$ (the dual space). For instance, in two dimensions it suffices to map a point $(\tilde{z},\tilde{t})$ to a line $L: \tilde{z}z-\tilde{t}$  and vice versa. Hence, the worst-case complexity of the $k$-levels is within a constant factor of the worst-case number of $k$-sets in point sets (\cite{dey1998,yu2012}).
 
 In  \cite{edelsbrunner1986} it is proved that, if the intersections of lines are computed dynamically in their manner, then Algorithm \ref{alg:Rider} constructs the boundary of the $k$-belt of $\mathcal{L}$ in $O(n_k\log^2 |\mathcal{L}|)$ time and $O(n_k+|\mathcal{L}|)$ space, where $n_k$ denotes the maximum number of edges bounding the $k$-belt of any set $\mathcal{L}$ with $|\mathcal{L}|$ lines. Upper bounds on $n_k$ are given in \cite{edelsbrunner1985} by transforming the geometric problem into  a combinatorial one. In particular, the upper bound established depends on $k$ and $|\mathcal{L}|$: $n_k=O(|\mathcal{L}|k^{\frac{1}{2}})$. \cite{dey1998,toth2000,nivasch2008} have subsequently improved the upper bound of the complexity on the $k$-set problem in the two-dimensional case. As for the development of efficient algorithms to compute $k$-envelopes or $k$-sets in higher dimensional spaces, they constitute a major challenge mainly because the computational complexity of $k$-envelopes increases exponentially with the increase in dimensionality (\cite{agarwal1995}). However, there exist algorithms for the cases $q=3$ (\cite{sharir2000, ramos1999, toth2000}),  $q=4$ (\cite{sharir2011}) and $q \ge 5$ (\cite{agarwal1998,halperin2017}).

As for the convexification procedure, the convex hull problem is a classical problem with immense literature.  The generalization of the Jarvis March to the multidimensional case is known as the \emph{graph transversal method} (and also as the \emph{gift-wrapping} algorithm), whereas the Graham Scan extension can be found under the name of \emph{incremental method}. For small dimensions ($q=2,3$), we encounter parallel algorithms in the literature that use the \emph{divide-and-conquer} paradigm to obtain better running times than the ones proposed here. However, we chose sequential ones for the ease of computation and explanation, and mainly because, given that the sets studied in our computational experiments are quite small, the difference in performance in the preprocessing steps is negligible. Surveys on algorithms and on the complexity of the convex hull problem can be found, for instance, in \cite{toth2017,preparata2012}. 

\end{document}